\documentclass[12pt, reqno]{amsart}
\usepackage{paralist}
\IfFileExists{ifsym.sty}{\usepackage[misc]{ifsym}}{} 
\usepackage{epsfig} 
\usepackage{epstopdf} 
\usepackage[colorlinks=true]{hyperref}
\usepackage{microtype}
\usepackage{amssymb,amsthm,amsmath}
\usepackage[numbers,sort&compress]{natbib}
\usepackage{amssymb,amsmath}
\usepackage{amsfonts}
\usepackage{mathrsfs}
\usepackage{latexsym}
\usepackage{amssymb}
\usepackage{amsthm}
\usepackage{pdfsync}
\usepackage{indentfirst}
\usepackage{xcolor}
\DeclareRobustCommand{\NSrev}[1]{\textcolor{blue}{#1}}

\usepackage{amsmath}
\usepackage{amsthm}

\newtheorem{theorem}{Theorem}[section]
\newcommand{\R}{\mathbb R}
\newcommand{\Y}{\mathcal Y}

\newtheorem{remark}[theorem]{Remark}
\newtheorem{proposition}[theorem]{Proposition}
\newtheorem{lemma}[theorem]{Lemma}
\newtheorem{definition}[theorem]{Definition}
\newtheorem{corollary}[theorem]{Corollary}

\newcommand{\beq}{\begin{equation}}
\newcommand{\eeq}{\end{equation}}
\newcommand{\ben}{\begin{eqnarray}}
\newcommand{\een}{\end{eqnarray}}

\numberwithin{equation}{section}

\title[Single-scale regularity for stationary NS]
{Single-scale regularity from Dirichlet energy for the
five-dimensional stationary Navier--Stokes equations}

\allowdisplaybreaks

\begin{document}

\author{Shuai~Li$^{1}$}
\address[Shuai~Li]{School of Mathematics and Statistics, Shandong University of Technology, Zibo, 255049,  China}
\email{leeshy@sdut.edu.cn}

\author{Jianbo~Yu$^{1,*}$}
\address[Jianbo~Yu]{School of Mathematical Sciences, Dalian University of Technology, Dalian, 116024,  China}
\email{yujb@mail.dlut.edu.cn}
\thanks{$^*$Author to whom correspondence should be addressed.}

	\date{\today}

\begin{abstract}
We establish an interior $\varepsilon$-regularity criterion for
suitable weak solutions of the unforced stationary Navier--Stokes
equations in an open set $\Omega\subset\mathbb R^5$.
There exists a universal constant $\varepsilon_\ast>0$ such that
$x_0$ is regular if, for some ball $B_R(x_0)\subset\Omega$ with $R>0$,
\[
\frac1R\int_{B_R(x_0)}|\nabla u|^2\,dx
\leq \varepsilon_\ast^2.
\]
The main difficulty is that the Dirichlet energy controls the velocity oscillation while its mean may be arbitrarily large. 
Subtracting the mean introduces a large constant drift into the perturbation equation.
Moreover, smallness in $H^1$ alone is insufficient to close a direct nonlinear elliptic estimate in five dimensions.
We overcome these difficulties through uniform Oseen estimates and a relative-energy comparison yielding quadratic improvement of the perturbation.
As a consequence, the singular set has zero one-dimensional upper Minkowski content on every compact subset of $\Omega$, and hence local upper Minkowski dimension at most one.
\end{abstract}

\maketitle
\thispagestyle{empty}
\medskip
\noindent\textbf{Keywords.}
Stationary Navier--Stokes equations,
single-scale regularity,
Dirichlet energy,
Oseen estimates,
Minkowski dimension.

\tableofcontents

\section{Introduction}
\subsection{Background and Main Results}

Consider the following system in $\Omega \subset \mathbb R^5$:
\begin{equation}\label{NS}
-\Delta u+(u\cdot\nabla)u+\nabla p=0,\qquad \nabla \cdot u=0,
\end{equation} 
where \(u:\Omega\to\mathbb{R}^{5}\) is the velocity field and \(p:\Omega\to\mathbb{R}\) is the pressure.

Recall that these so-called $\varepsilon$-regularity criteria can be traced back to the well-known work by Caffarelli-Kohn-Nirenberg \cite{CKN82} for the analysis of suitable weak solutions of the three dimensional time-dependent Navier-Stokes equations, where they showed that the set $S$ of possible interior singular points of a suitable weak solution has zero one-dimensional parabolic Hausdorff measure by improving Scheffer's results in \cite{Sch76,Sch77,Sch80}.  A simplified proof was given later by Lin in \cite{Lin98} with $p \in L^\frac32(Q_1)$. In \cite{TianXin}, Tian and Xin proved
a local uniform $L^\infty$ bound of $\nabla u$ for smooth solutions
when either the scaled local $L^2$ norm of the vorticity
or the scaled local total energy is small. In \cite{Vasseur}, Vasseur proved that, for \(\tilde q>1\), if
\begin{align*}
\sup_{t\in[-1,0]}
\int_{B_1}|u(t,x)|^2\,dx
+\int_{Q_1}|\nabla u|^2\,dx\,dt
+\left[\int_{-1}^{0}\|p(t,\cdot)\|_{L^1(B_1)}^{\tilde q}\,dt\right]^\frac{2}{\tilde{q}}
\end{align*}
is sufficiently small, then \(u\) is regular in \(Q_{\frac12}\),
where \(B_r\subset\mathbb R^3\) and \(Q_r=B_r\times(-r^2,0)\). Wang and Zhang \cite{WangZhang14} proved regularity in the backward parabolic cylinder $Q_{\frac12}$ when the following quantity is sufficiently small:
\begin{align*}
\sup_{t\in[-1,0]} \int_{B_1} |u(t, x)|^2 \, dx + \int_{-1}^0 \left(\int_{B_1} |u|^4 \, dx\right)^\frac12 \, dt + \int_{-1}^0\left(\int_{B_1} |p|^2 \, dx\right)^\frac12dt.
\end{align*}
For more regularity results on the time-dependent Navier-Stokes equations, we refer the reader to \cite{HWZ19,DongWang20,Kukavica,DongDu07,SeShSo} and the references therein.
Wolf \cite{Wolf10,Wolf15} proved that the single-scale condition
\begin{align*}
\int_{Q_1} |u|^3 \, dxdt \leq \varepsilon_{Wolf}
\end{align*}
implies local regularity, using a local pressure decomposition.
See also interesting results in \cite{Kwon23,ChaeWolf17ARMA,LiWangZhou23}.

In this paper, however, we focus on the stationary Navier-Stokes system. Although a stationary solution can be viewed formally as a time-independent
solution of the evolution equations, the stationary problem has its own analytic features, and the arguments for the time-dependent system do not directly yield the result needed below.
Struwe established $\varepsilon$-regularity for the five-dimensional stationary equations in \cite{St88} using methods from the regularity theory of elliptic systems (cf. Morrey \cite{Morrey66} and Giaquinta \cite{G83}). He proved that $0$ is a regular point if 
\begin{align*}
\limsup_{r \to 0} r^{-1} \int_{B_r} |\nabla u|^2 dx \leq \varepsilon_{St}.
\end{align*}
Later, the result of Struwe was extended to the boundary case by Kang \cite{Kang04}. 
In \cite{JiaSv18}, under the assumption that $\|u\|_{L^3(B_1)}$ is sufficiently small, Jia and \v{S}ver\'ak proved that the suitable weak solution to \eqref{NS} is smooth in \(B_{\frac12}\), both in \(\mathbb{R}^{5}\) and \(\mathbb{R}^{6}\).
Cui \cite{Cui2023} extended the single-scale $L^3$ criterion
of Jia and \v{S}ver\'ak \cite{JiaSv18} to $L^q$ criteria
with $\frac52<q<\frac{10}{3}$ in five dimensions, without requiring
pressure smallness.
In the unforced case, the smallness of the velocity alone
is sufficient to conclude that $0$ is a regular point;
in particular,
\begin{align*}
\int_{B_1} |u|^3 dx <\varepsilon_{\mathrm{cui}} \Rightarrow u \in L^\infty(B_\frac18).
\end{align*}
The boundary case is also treated in \cite{Cui2023}. 
A closely related result was obtained by Li and Wang \cite{LiWang23} for the six-dimensional stationary equations.
In the unforced interior case, their single-scale criterion requires
\[
\left[ 1+\left(r^{-5}\int_{B_r}|u|\,dx\right)^2 \right] r^{-2}\int_{B_r}|\nabla u|^2\,dx <\varepsilon_{LW}
\]
at some radius $r>0$.
Thus, this criterion retains a factor involving the velocity.
For the boundary problem with zero velocity trace, they also obtained a single-scale criterion involving only the Dirichlet energy in the unforced case. 
For more relevant research about the higher dimensional stationary Navier-Stokes equations, we refer readers to \cite{DongGu14,DongStrain12,FrRu94,FrRu95,JiaSv18,LuoX19,LiuWang18} and the references therein.

The question discussed here concerns the five-dimensional interior problem: {\bf Does sufficiently small scale-invariant Dirichlet energy on a single ball imply regularity, with a universal threshold and without any bound on the mean velocity?}
Theorem~\ref{thm:main} gives an affirmative answer.


We use the following definition of an interior suitable weak solution.
\begin{definition} \label{Def:SWS}
Let  $\Omega\subset\mathbb{R}^5$ be an open domain. The pair $(u,p)$ is said to be a suitable weak solution to the steady Navier--Stokes equations \eqref{NS} in $\Omega$ if the following conditions hold.

(i)\,$u\in H^1(\Omega),\,p\in L^{\frac 32}(\Omega),\,\nabla p \in L^{\frac 65}(\Omega)$;

(ii)\,$(u,p)$ satisfies the equations \eqref{NS} in the sense of distributions;

(iii)\,$u$ and $p$ satisfy the local energy inequality
\ben\label{eq:local energy}
2\int_{\Omega}|\nabla u|^2\phi dx \leq\int_{\Omega}\big[|u|^2\triangle\phi+u\cdot\nabla\phi(|u|^2+2p)\big] dx
\een
for every nonnegative $\phi\in C_c^\infty(\Omega)$.
\end{definition}
A point is regular if $u$ is essentially bounded in a neighborhood of it. In the unforced stationary system this implies interior smoothness by the Stokes estimates.
Our main result is as follows:

\begin{theorem}\label{thm:main}
Let \((u,p)\) be a suitable weak solution of \eqref{NS} in \(\Omega\).
There exists a universal constant \(\varepsilon_*>0\) such that,
if \(B_R(x_0)\subset\Omega\) for some \(R>0\) and
\begin{equation}\label{oneScale}
\frac1R\int_{B_R(x_0)}|\nabla u|^2\,dx
\leq \varepsilon_*^2,
\end{equation}
then \(x_0\) is a regular point.
\end{theorem}

\begin{corollary} \label{minkowski} 
Let \((u,p)\) be a suitable weak solution of $\eqref{NS}$ in \(\Omega\). Let \(\mathcal S\) be the set of singular points.
For every compact set \(K\Subset\Omega\), set
\begin{align*}
S_K=\mathcal S\cap K, \qquad (S_K)_r =\{x\in\mathbb R^5:\operatorname{dist}(x,S_K)<r\}.
\end{align*}
Then
\begin{equation}\label{MinkowskiConclusion}
\mathcal M^1(S_K)=0, \qquad \overline{\dim}_{\mathrm M}S_K\le1,
\end{equation}
where the one-dimensional upper Minkowski content is defined by
\begin{align*}
\mathcal M^1(S_K) :=\limsup_{r\to0} \frac{|(S_K)_r|}{\omega_4r^4}, \qquad \omega_4=|B_1^4|.
\end{align*}
More precisely, for all sufficiently small \(r>0\),
\begin{equation}\label{tubeEstimate}
|(S_K)_r| \le\frac{C}{\varepsilon_\ast^2}\,r^4 \int_{(S_K)_r}|\nabla u|^2\,\mathrm dx,
\end{equation}
where \(C\) depends only on the dimension and \(\varepsilon_\ast\) is the constant in Theorem \ref{thm:main}.
\end{corollary}

\begin{remark}
For each compact set \(K\Subset\Omega\), estimate \eqref{tubeEstimate} provides a quantitative bound on the volume of the neighborhood of \(S_K\):
\begin{align*}
\frac{|(S_K)_r|}{\omega_4 r^4} \leq \frac{C}{\varepsilon_*^2} \int_{(S_K)_r} |\nabla u|^2\,dx
\qquad\text{for all sufficiently small }r>0.
\end{align*}
Since \(S_K\) is compact and has zero five-dimensional Lebesgue measure, the neighborhoods \((S_K)_r\) decrease to \(S_K\) as \(r\to0\). The local integrability of \(|\nabla u|^2\) therefore implies
\begin{align*}
\lim_{r\to 0} \int_{(S_K)_r}|\nabla u|^2\,dx=0.
\end{align*}
Consequently,
\begin{align*}
\mathcal{M}^{1}(S_K) = \limsup_{r\to 0} \frac{|(S_K)_r|}{\omega_4 r^4} =0.
\end{align*}
In particular, \(\mathcal{H}^{1}(S_K)=0\), and hence \(\mathcal{H}^{s}(S_K)=0\) for every \(s\geq1\). 
Rather, it shows that the possible singular set has zero one-dimensional Hausdorff measure, while \eqref{tubeEstimate} quantifies the decay of the volume of its tubular neighborhoods in terms of the Dirichlet energy near \(S_K\).
\end{remark}

\subsection{Motivation and Outline of Proof}
Classical iteration arguments (see \cite{St88}) require the scale-invariant energy to be small at every scale. 
Smallness at one fixed scale does not directly provide this
control, since the normalized energy may increase at smaller scales.
Thus, in order to establish a regularity criterion based solely on the scale-invariant Dirichlet energy on one ball, deriving regularity from smallness on a single ball requires an additional mechanism to propagate the information at the initial scale to the interior scales. 
A main methodological contribution is the local quadratic improvement established in Proposition~\ref{improvement}.
The proof uses a relative-energy argument to compare the suitable weak solution with a smooth comparison field.
The initial energy controls the velocity oscillation, but not its mean. 
After subtracting the mean, a small $H^1$ error remains, together with a possibly large constant drift. 
We retain this drift in the principal linear operator and prove a uniform estimate independent of its size.
The nonlinear difficulty is that smallness in \(H^1\) is insufficient to close a direct elliptic estimate in five dimensions.

 At each stage of the construction, we write the solution in the form $u = a e_1 + U + w,$ where $a e_1 + U$ is a smooth exact solution of the stationary Navier-Stokes equations and $w$ is the weak error. 
The improved decomposition is then rescaled to the unit cylinder and the procedure is repeated. Schematically, if $\varepsilon_j$ measures the weak error, $\delta_j$ measures the size of $U_j$ in the mixed norm introduced below, and $a_j$ the constant drift at the $j$-th stage, the estimates take the form
\[
    \varepsilon_{j+1} \leq C_{\mathrm{imp}}\theta^{-\frac32}\varepsilon_j^2,
    \qquad
    \delta_{j+1} \leq \delta_j + C_{\mathrm{imp}}\varepsilon_j,
    \qquad
    a_{j+1} = \theta a_j,
\]
for a fixed $\theta \in (0,1)$. The quadratic decay of the error makes the increments in the background size summable, while rescaling reduces the drift geometrically. Consequently, sufficiently small initial Dirichlet energy guarantees that every required correction is admissible, independently of the number of steps.

After finitely many steps, the rescaled drift is bounded. 
At that scale, the velocity is the sum of a bounded constant field, a small smooth background, and a small weak error. 
This yields, on a smaller ball, the velocity smallness needed to apply a known one-scale velocity criterion, which implies regularity near the center of the original ball.

The rest of the paper is organized as follows. Section~2 introduces the scaling, mixed norms, and analytic tools used throughout the paper. Section~3 is devoted to establishing full-space and interior estimates for the Oseen system, together with estimates for its linearization around a small smooth background, uniformly with respect to the large constant drift. Section~4 constructs smooth nonlinear corrections and derives a relative energy inequality, leading to the quadratic improvement of the weak perturbation.
In Section~5 we prove the single-scale regularity criterion by a finite iteration and establish the Minkowski dimension and content estimates for the singular set. Some auxiliary estimates and the divergence-free extension construction are added in the Appendix.

\section{Preliminaries}

Set $B_R(x_0) = \{x \in \mathbb R^5; |x - x_0| < R\}$ and $B_R = \{x \in \mathbb R^5; |x| < R\}$. 
By translation, take $x_0=0$ in Theorem \ref{thm:main}.
Define a scaling transformation
\begin{align*}
v(y) = R u(R y), \quad q(y) = R^2 p(Ry), \quad y \in B_1\NSrev{.}
\end{align*}
If $(u,p)$ solves \eqref{NS}, then $(v,q)$ also solves \eqref{NS}. The local energy inequality is preserved under this scaling, and
\begin{align*}
\int_{B_1} |\nabla v|^2 dy = R^{-1} \int_{B_R} |\nabla u|^2 dx.
\end{align*}

Next, we write \(x=(s,y)\in\mathbb R\times\mathbb R^4\),  and let \(B_r^4=\{x\in\mathbb R^4:|x|<r\}\), and set
\begin{align*}
Q_r = (-r,r)\times B_r^4,\qquad L_a= -\Delta + a \partial_s.
\end{align*}
Here \(s\) is a spatial coordinate. For \(a>0\), define
\begin{equation}\label{Ydef}
\Y_a(f;Q_r)= \|f\|_{L^\infty((-r,r);L^4(B_r^4))} +\|\nabla f\|_{L^\infty((-r,r);L^2(B_r^4))} +a^{-\frac12}\|\nabla^2f\|_{L^2(Q_r)}.
\end{equation}
The full-space version is
\begin{align*}
\Y_a(f)=\|f\|_{L^\infty(\mathbb R;L^4(\mathbb R^4))} +\|\nabla f\|_{L^\infty(\mathbb R;L^2(\mathbb R^4))} +a^{-\frac12}\|\nabla^2f\|_{L^2(\mathbb R^5)}.
\end{align*}
For $\theta \in (0,1)$, under the transform $f_\theta(x) = \theta f(\theta x)$, there holds
\begin{align*}
\|f_\theta\|_{L^\infty((-1,1);L^4(B_1^4))} &=\|f\|_{L^\infty((-\theta,\theta);L^4(B_\theta^4))},\\
\|\nabla f_\theta\|_{L^\infty((-1,1);L^2(B_1^4))} &=\|\nabla f\|_{L^\infty((-\theta,\theta);L^2(B_\theta^4))},\\
\|\nabla^2f_\theta\|_{L^2(Q_1)} &=\theta^{\frac12}\|D^2f\|_{L^2(Q_\theta)}.
\end{align*}
Moreover,
\begin{align*}
\|f_\theta\|_{H^1(Q_1)}^2 =\theta^{-3}\|f\|_{L^2(Q_\theta)}^2 +\theta^{-1}\|\nabla f\|_{L^2(Q_\theta)}^2 \leq \theta^{-3}\|f\|_{H^1(Q_\theta)}^2.
\end{align*}
Hence
\begin{equation}\label{Yscale}
\Y_{\theta a}(f_\theta;Q_1)=\Y_a(f;Q_\theta), \qquad \|f_\theta\|_{H^1(Q_1)} \le\theta^{-\frac32}\|f\|_{H^1(Q_\theta)}.
\end{equation}

Let us introduce the following lemma. 

\begin{lemma}[Hardy-Littlewood-Sobolev inequality, \cite{Stein1970}]\label{lem:HLSine}
Let $0 < \alpha < n$ and let $1 < p < \frac n \alpha$. Define
the fractional integral operator
\[
I_{\alpha}f(x)
=
c_{n,\alpha}
\int_{\mathbb{R}^{n}}
\frac{f(y)}{|x-y|^{n-\alpha}}\,dy,
\]
where $c_{n,\alpha}$ is a dimensional normalization constant. If
\[ \frac{1}{q} = \frac{1}{p} - \frac{\alpha}{n}, \]
then there exists a constant
$C=C(n,\alpha,p)>0$ such that
\[ \|I_{\alpha}f\|_{L^{q}(\mathbb{R}^{n})} \leq C\|f\|_{L^{p}(\mathbb{R}^{n})}.\]
\end{lemma}

\begin{lemma}[Iteration lemma; cf.~\cite{G83}, Ch.~V, Lemma~3.1]
\label{lem:Gia}
Let $f(r)$ be a non-negative bounded function on $[R_0,R_1] \subset \R_+$. If there are nonnegative constants $A,B,D$ and positive exponents $b<a$ and a parameter $\theta \in (0,1)$ such that for all $R_0 \le \rho < \tau \le R_1$
\[
f(\rho) \le \theta f(\tau) + \frac{A}{(\tau-\rho)^{a}} + \frac{B}{(\tau-\rho)^{b}} + D,
\]
then for all $R_0 \le \rho < \tau \le R_1$
\[
f(\rho) \leq C(a,\theta) \left[\frac{A}{(\tau-\rho)^{a}} + \frac{B}{(\tau-\rho)^{b}} + D\right].
\]
\end{lemma}

\begin{lemma}[Marcinkiewicz multiplier theorem, \cite{GS17}]\label{lem:mmt}
Let $m:\mathbb{R}^{d}\to\mathbb{C}$ be a bounded function which is
smooth on each open coordinate orthant. For
$J\subset\{1,\ldots,d\}$, write
\[
\partial_J=\prod_{j\in J}\partial_{\xi_j},
\qquad
\partial_{\varnothing}m=m.
\]
Assume that
\[
M_m:=
\max_{J\subset\{1,\ldots,d\}}
\sup_{\xi_1\cdots\xi_d\neq 0}
\left|
\left(\prod_{j\in J}\xi_j\right)
\partial_Jm(\xi)
\right|
<\infty.
\]
Then the Fourier multiplier operator
\[
T_mf
=
\mathcal{F}^{-1}\bigl(m\,\widehat f\bigr)
\]
extends boundedly on $L^p(\mathbb{R}^d)$ for every
$1<p<\infty$, and
\[
\|T_mf\|_{L^p(\mathbb{R}^d)}
\leq
C_{d,p}M_m\|f\|_{L^p(\mathbb{R}^d)}.
\]
The values of $m$ on the coordinate hyperplanes are irrelevant.
\end{lemma}

\section{Uniform estimates for the Oseen system} 
\subsection{Full-space estimates}

We establish full-space estimates for the Oseen resolvent,
uniformly for \(a\ge1\), and prove the invertibility of
\(I+B_U\) under a smallness condition on the background \(U\).

The unitary Fourier transform is as follows:
\begin{align*}
\widehat f(\xi)=(2\pi)^{-\frac52}\int_{\mathbb R^5}e^{-ix\cdot\xi}f(x)\,\mathrm dx,
\qquad
\mathcal F_yf(s,\xi')=(2\pi)^{-2}\int_{\mathbb R^4} e^{-iy\cdot\xi'}f(s,y)\,\mathrm dy.
\end{align*}
Let $\mathbb P$ denote the full-space Leray projection, with matrix symbol
$\mathsf P(\xi)=I-\xi\otimes\xi/|\xi|^2$ for $\xi\ne0$. In particular,
$|\mathsf P(\xi)\eta|\le|\eta|$. Set
\begin{align*}
\mathcal R_a=L_a^{-1}\mathbb P,\qquad S_a T=\mathcal R_a \nabla \cdot T,
\qquad \xi=(\tau,\xi'),\quad k=|\xi'|,\quad
D_a=\tau^2+k^2+ia\tau.
\end{align*}
Thus \(\widehat{\mathcal R_a F}=D_a^{-1}\mathsf P(\xi)\widehat F\) and
\begin{align*}
\widehat{ S_a T}_i(\xi)
=D_a^{-1}\sum_{{\ell},j=1}^5\mathsf P_{i{\ell}}(\xi)i\xi_j\widehat T_{{\ell} j}(\xi).
\end{align*}

\begin{lemma}\label{resolvent}
For \(a\ge1\), the following estimates have constants independent of \(a\):
\begin{align}
\|\nabla S_a T\|_{L^p(\mathbb R^5)} &\le C_p\|T\|_{L^p(\mathbb R^5)},\qquad 1<p<\infty, \label{CZ}\\*
\sqrt a\,\| S_a T\|_{L^2(\mathbb R^5)} &\le C\|T\|_{L^1(\mathbb R;L^2(\mathbb R^4))}, \label{driftL2}\\*
\Y_a(\mathcal R_a F)&\le Ca^{-\frac12}\|F\|_{L^2(\mathbb R^5)}. \label{RF}
\end{align}
If \(T\in L^2(\mathbb R^5)\) is supported in a fixed bounded measurable set \(K\), then
\begin{equation}\label{compactT}
\| S_a T\|_{H^1(\mathbb R^5)}
\le C(1+|K|^{\frac15})\|T\|_{L^2(\mathbb R^5)}.
\end{equation}
\end{lemma}

\begin{proof}

Put \(A_\xi=|\xi|^2\) and \(\xi=(\tau,\xi')=(\xi_1,...,\xi_5)\), so that \(D_a=A_\xi+ia\xi_1\), and consider \(m_a=A_\xi D_a^{-1}\). For a subset \(J\subset\{1,\ldots,5\}\), write \(\partial_J=\prod_{j\in J}\partial_{\xi_j}\), with \(\partial_\varnothing\) the identity. 
We now check every required derivative of \(m_a\). 

{\underline{\it proof of \eqref{CZ}.}}
Since
\begin{align*}
|D_a|=(A_\xi^2+a^2\xi_1^2)^{\frac12}\ge A_\xi,\qquad
|D_a|\ge a|\xi_1|,
\end{align*}
we have
\begin{align*}
\left|\frac{A_\xi}{D_a}\right|\le1,\qquad
\left|\frac{\xi_j\partial_{\xi_j}D_a}{D_a}\right|
=\left|\frac{2\xi_j^2+ia\xi_1\mathbf1_{\{j=1\}}}{D_a}\right|\le3,
\qquad
\left|\frac{\xi_j\partial_{\xi_j}A_\xi}{D_a}\right|
=\left|\frac{2\xi_j^2}{D_a}\right|\le2.
\end{align*}
Here \(\partial_{\xi_j}D_a=2\xi_j+ia\mathbf1_{\{j=1\}}\), while
\(\partial_{\xi_i}\partial_{\xi_j}D_a=0\) when \(i\ne j\).
This absence of mixed derivatives gives the exact formula
\begin{equation}\label{inverseMixed}
\partial_J(D_a^{-1})
=(-1)^{|J|}|J|!D_a^{-|J|-1}
 \prod_{j\in J}\partial_{\xi_j}D_a.
\end{equation}
Formula \eqref{inverseMixed} follows by induction, since
\(\partial_{\xi_i}\partial_{\xi_j}D_a=0\) for \(i\ne j\).
Likewise, mixed derivatives of \(A_\xi\) of order at least two
vanish, so
\begin{align*}
\partial_Jm_a
=A_\xi\partial_J(D_a^{-1})
+\sum_{{\ell}\in J}
(\partial_{\xi_{\ell}}A_\xi)
\partial_{J\setminus\{{\ell}\}}(D_a^{-1}).
\end{align*}
For \(n=|J|\ge1\), the preceding bounds therefore give
\begin{align*}
\left|
\left(\prod_{j\in J}\xi_j\right)\partial_Jm_a
\right|
\le n!3^n+2n(n-1)!3^{n-1}.
\end{align*}
Together with \(|m_a|\le1\), this verifies the Marcinkiewicz
conditions uniformly in \(a\), and hence the uniform
\(L^p(\mathbb R^5)\) multiplier bound.
Finally, the symbol mapping \(T_{{\ell} j}\) to \(\partial_m( S_a T)_i\) is
\begin{align*}
-\frac{\xi_m\xi_j}{D_a}\mathsf P_{i{\ell}}(\xi)
=\left(-\frac{\xi_m\xi_j}{|\xi|^2}\right)m_a(\xi)\mathsf P_{i{\ell}}(\xi).
\end{align*}
By Lemma \ref{lem:mmt}, \eqref{CZ} holds.

{\underline{\it proof of \eqref{driftL2}.}}
The Fourier transform in the longitudinal variable satisfies
\begin{align*}
|\widehat T(\tau,\xi')|
=\left|(2\pi)^{-\frac12}\int_{\mathbb R}
 e^{-is\tau}\mathcal F_yT(s,\xi')\,\mathrm ds\right|
\le(2\pi)^{-\frac12}\int_{\mathbb R}|\mathcal F_yT(s,\xi')|\,\mathrm ds.
\end{align*}
For \(k>0\), set
\(\mu_+=(\sqrt{a^2+4k^2}+a)/2\) and
\(\mu_-=(\sqrt{a^2+4k^2}-a)/2\). Both are positive, and
\begin{align*}
\mu_+\mu_-=k^2,\qquad \mu_++\mu_-=\sqrt{a^2+4k^2},\qquad
\mu_+^2+\mu_-^2=a^2+2k^2.
\end{align*}
Factoring the denominator and integrating the partial fractions in one
calculation gives
\begin{equation}\label{tauint}
\begin{aligned}
\int_{\mathbb R}\frac{\tau^2+k^2}{(\tau^2+k^2)^2+a^2\tau^2}\,\mathrm d\tau
&=\int_{\mathbb R}\frac{\tau^2+\mu_+\mu_-}
 { (\tau^2+\mu_+^2)(\tau^2+\mu_-^2)}\,\mathrm d\tau\\
&=\frac1{\mu_++\mu_-}\int_{\mathbb R}
 \left(\frac{\mu_+}{\tau^2+\mu_+^2}
       +\frac{\mu_-}{\tau^2+\mu_-^2}\right)\,\mathrm d\tau\\
&=\frac{\pi+\pi}{\mu_++\mu_-}
 =\frac{2\pi}{\sqrt{a^2+4k^2}}\le\frac{2\pi}{a}.
\end{aligned}
\end{equation}
At \(k=0\), cancellation of \(\tau^2\) away from \(\tau=0\) gives
\(\int_{\mathbb R}(\tau^2+a^2)^{-1}\,\mathrm d\tau=\frac{\pi}{a}\), so the same
upper bound holds. The value of the integrand at a single point is immaterial.
Using the symbol of \( S_a\), \(|\mathsf P(\xi)|\le1\), Plancherel's theorem,
and \eqref{tauint}, we obtain
\begin{align*}
\| S_a T\|_{L^2(\mathbb R^5)}^2
&\le\int_{\mathbb R^4}\int_{\mathbb R}
 \frac{\tau^2+|\xi'|^2}{|D_a|^2}|\widehat T(\tau,\xi')|^2
 \,\mathrm d\tau\,\mathrm d\xi' \le\frac C a\int_{\mathbb R^4}
 \left(\int_{\mathbb R}|\mathcal F_yT(s,\xi')|\,\mathrm ds\right)^2
 \,\mathrm d\xi'\\
&\le\frac C a\left[\int_{\mathbb R}
 \left(\int_{\mathbb R^4}|\mathcal F_yT(s,\xi')|^2
 \,\mathrm d\xi'\right)^{\frac12}\,\mathrm ds\right]^2\\
&=\frac C a\left[\int_{\mathbb R}
 \left(\int_{\mathbb R^4}|T(s,y)|^2\,\mathrm dy\right)^{\frac12}
 \,\mathrm ds\right]^2
 =\frac C a\|T\|_{L^1(\mathbb R;L^2(\mathbb R^4))}^2.
\end{align*}
This proves \eqref{driftL2}.

{\underline{\it proof of \eqref{RF}.}}
First take \(F\in C_c^\infty(\mathbb R^5)\) and let \(z=\mathcal R_a F\). Partial Fourier inversion gives
\begin{align*}
\mathcal F_y[\nabla z](s,\xi')
=(2\pi)^{-\frac12}\int_{\mathbb R}e^{is\tau}
 \frac{i\xi\otimes \mathsf P(\xi)\widehat F(\tau,\xi')}{D_a}\,\mathrm d\tau.
\end{align*}
Cauchy--Schwarz in \(\tau\), followed by \eqref{tauint}, yields
\begin{align*}
|\mathcal F_y[\nabla z](s,\xi')|^2
\le\frac1{2\pi}
 \left(\int_{\mathbb R}\frac{\tau^2+|\xi'|^2}{|D_a|^2}\,\mathrm d\tau\right)
 \left(\int_{\mathbb R}|\widehat F(\tau,\xi')|^2\,\mathrm d\tau\right)
\le\frac1a\int_{\mathbb R}|\widehat F(\tau,\xi')|^2\,\mathrm d\tau.
\end{align*}
The right-hand side is independent of \(s\). Integrating in \(\xi'\)
and using Plancherel in both sets of variables therefore gives
\begin{align*}
&\|\nabla z\|_{L^\infty(\mathbb R;L^2(\mathbb R^4))}^2
\le\frac1a\int_{\mathbb R^4}\int_{\mathbb R}
 |\widehat F(\tau,\xi')|^2\,\mathrm d\tau\,\mathrm d\xi'
 =\frac1a\|F\|_{L^2(\mathbb R^5)}^2,\\*
&\|z\|_{L^\infty(\mathbb R;L^4(\mathbb R^4))}
\le C\|\nabla_y z\|_{L^\infty(\mathbb R;L^2(\mathbb R^4))}
 \le Ca^{-\frac12}\|F\|_{L^2(\mathbb R^5)},\\*
&\|D^2z\|_{L^2(\mathbb R^5)}^2
\le\int_{\mathbb R^5}\frac{|\xi|^4}{|D_a|^2}|\widehat F(\xi)|^2
 \,\mathrm d\xi
 \le\int_{\mathbb R^5}|\widehat F(\xi)|^2\,\mathrm d\xi
 =\|F\|_{L^2(\mathbb R^5)}^2.
\end{align*}
Here, we use the relation \(|D_a|\ge|\xi|^2\) and the homogeneous Sobolev inequality on \(\mathbb{R}^4\).
Summing the three estimates and using density in
\(L^2(\mathbb R^5)\) yields \eqref{RF}.

{\underline{\it proof of \eqref{compactT}.}}
The pointwise Fourier-symbol bound \(|\xi|/|D_a|\le|\xi|^{-1}\), Plancherel, and Lemma \ref{lem:HLSine} give 
\begin{align*}
\| S_a T\|_{L^2(\mathbb R^5)}
&\le\left(\int_{\mathbb R^5}|\xi|^{-2}|\widehat T(\xi)|^2
 \,\mathrm d\xi\right)^{\frac12}
 =\|(-\Delta)^{-\frac12}T\|_{L^2(\mathbb R^5)}\\
&\le C\|T\|_{L^{\frac{10}{7}}(\mathbb R^5)}
 \le C|K|^{\frac15}\|T\|_{L^2(\mathbb R^5)}.
\end{align*}
Combining this with \eqref{CZ} for \(p=2\) proves \eqref{compactT}.
\end{proof}

\begin{corollary}
\label{cor:bilinear-estimate}
Let $a\geq 1$, and let $V,W$ be smooth vector fields on
$\mathbb{R}^{5}$ with $\Y_a(V)+\Y_a(W)<\infty$. Then
\begin{align*}
\Y_a\bigl(\mathcal R_a(V\cdot\nabla W)\bigr) &\leq C a^{-\frac12} \|V\cdot\nabla W\|_{L^2(\mathbb{R}^5)} \\
&\leq C a^{-\frac12} \|V\|_{L^\infty(\mathbb{R};L^4(\mathbb{R}^4))}  \|D^2W\|_{L^2(\mathbb{R}^5)}
\leq C \Y_a(V)\Y_a(W),
\end{align*}
where $C$ is independent of $a$.

For the same fields,
\[
\Y_a\bigl(\mathcal R_a \nabla \cdot (V\otimes W)\bigr) \leq C \Y_a(V)\Y_a(W),
\]
without any divergence-free assumption on $V$ or $W$.
\end{corollary}

\begin{lemma}\label{cross}
Let \(a\ge1\) and \(U\in L^\infty(\mathbb R;L^4(\mathbb R^4))\). Set
\begin{align*}
\beta=\|U\|_{L^\infty(\mathbb R;L^4(\mathbb R^4))},
\qquad
B_Uz=S_a(U\otimes z+z\otimes U).
\end{align*}
Then
\begin{align}
\|B_Uz\|_{L^2(\mathbb R^5)} &\le C\beta\|z\|_{L^2(\mathbb R^5)}, \label{cross2}\\
\|\nabla B_Uz\|_{L^{\frac53}(\mathbb R^5)} &\le C\beta\|\nabla z\|_{L^{\frac53}(\mathbb R^5)}, \label{cross53}
\end{align}
where \(C\) is independent of \(a\) and \(U\).
Moreover, if \(C\beta\le\frac12\), then \(I+B_U\) is invertible on \(L^2(\mathbb R^5)\cap\dot W^{1,\frac53}(\mathbb R^5).\)
More precisely, for every \(z_0\) in this intersection, the unique solution of \((I+B_U)z=z_0\) satisfies
\begin{align*}
\begin{aligned}
\sqrt a\,\|z\|_{L^2(\mathbb R^5)} +\|\nabla z\|_{L^{\frac53}(\mathbb R^5)} \le 2\Bigl( \sqrt a\,\|z_0\|_{L^2(\mathbb R^5)} +\|\nabla z_0\|_{L^{\frac53}(\mathbb R^5)} \Bigr).
\end{aligned}
\end{align*}
\end{lemma}

\begin{proof}
We first establish the two operator estimates and then construct the inverse. 
Write \(\xi=(\tau,\xi')\) and
\(D_a=|\xi|^2+ia\tau\). Since the Leray symbol has operator norm
at most one, the symbol of \(S_a\) has size at most
\begin{align*}
\frac{|\xi|}{|D_a|}
=\frac{|\xi|}{\sqrt{|\xi|^4+a^2\tau^2}}
\le\frac1{|\xi|}
\le\frac1{|\xi'|}
\qquad(\xi'\ne0).
\end{align*}
Consequently, for a smooth compactly supported tensor \(T\), Plancherel's theorem and Lemma \ref{lem:HLSine} give
\begin{align*}
\|S_a T\|_{L^2(\mathbb R^5)}^2
\le C\int_{\mathbb R^4}\int_{\mathbb R}
\frac{|\widehat T(\tau,\xi')|^2}{|\xi'|^2}
\,\mathrm d\tau\,\mathrm d\xi'
&=C\int_{\mathbb R}
\|I_1^yT(s,\cdot)\|_{L^2(\mathbb R^4)}^2
\,\mathrm ds\\
&\le C\int_{\mathbb R}
\|T(s,\cdot)\|_{L^{\frac43}(\mathbb R^4)}^2
\,\mathrm ds.
\end{align*}
Taking \(T=U\otimes z+z\otimes U\), H\"older's inequality on transverse slices yields
\begin{align*}
\|B_Uz\|_{L^2(\mathbb R^5)}
&\le C
\left[
\int_{\mathbb R}
\|U\otimes z+z\otimes U\|_{L^{\frac43}(\mathbb R^4)}^2
\,\mathrm ds
\right]^{\frac12}\\
&\le C
\left[
\int_{\mathbb R}
\|U(s,\cdot)\|_{L^4(\mathbb R^4)}^2
\|z(s,\cdot)\|_{L^2(\mathbb R^4)}^2
\,\mathrm ds
\right]^{\frac12}\\
&\le C\beta
\left[
\int_{\mathbb R}\int_{\mathbb R^4}
|z(s,y)|^2\,\mathrm dy\,\mathrm ds
\right]^{\frac12}
=C\beta\|z\|_{L^2(\mathbb R^5)}.
\end{align*}
This proves \eqref{cross2}.

By \eqref{CZ}, H\"older's inequality and the embedding $\dot W^{1,\frac53}(\mathbb R^4)\hookrightarrow L^{\frac{20}{7}}(\mathbb R^4)$, we have
\begin{align*}
\|\nabla B_Uz\|_{L^{\frac53}(\mathbb R^5)}
&\le C
\|U\otimes z+z\otimes U\|_{L^{\frac53}(\mathbb R^5)}\\
&\le C\beta
\left[
\int_{\mathbb R}\int_{\mathbb R^4}
|\nabla_yz(s,y)|^{\frac53}
\,\mathrm dy\,\mathrm ds
\right]^{\frac35}\le C\beta\|\nabla z\|_{L^{\frac53}(\mathbb R^5)}.
\end{align*}
This proves \eqref{cross53}.

Set \(X_a=L^2(\mathbb R^5)\cap\dot W^{1,\frac53}(\mathbb R^5)\) with the norm
\begin{align*}
\|z\|_{X_a} = \sqrt a\,\|z\|_{L^2(\mathbb R^5)} +\|\nabla z\|_{L^{\frac53}(\mathbb R^5)}.
\end{align*}
Combining \eqref{cross2} and \eqref{cross53} gives
\begin{align*}
\|B_Uz\|_{X_a}\le C\beta\|z\|_{X_a}.
\end{align*}
If \(C\beta\le\frac12\), the Neumann series converges in operator
norm on \(X_a\), with
\begin{align*}
(I+B_U)^{-1}=\sum_{m=0}^{\infty}(-B_U)^m,
\qquad
\|(I+B_U)^{-1}\|_{X_a\to X_a}
\le\frac1{1-C\beta}\le2.
\end{align*}
Thus \(z=(I+B_U)^{-1}z_0\) is the unique solution in \(X_a\).
Expanding the \(X_a\) norm gives the stated estimate, while
\begin{align*}
\|z\|_{L^{\frac52}(\mathbb R^5)}
\le C\|\nabla z\|_{L^{\frac53}(\mathbb R^5)}
\end{align*}
follows from the homogeneous Sobolev inequality.
\end{proof}

\subsection{Interior estimates for homogeneous and linearized equations}

We prove interior \(\Y_a\) estimates for the homogeneous Oseen
system and its linearization around a small smooth background,
with constants uniform for \(a\ge1\).

\begin{lemma}\label{homogeneous}
Let \(0<r<R\le1\), \(a\ge1\), and suppose that \(h\in H^1(Q_R)\) and a distribution \(\pi\) satisfy
\begin{align}\label{eq:h}
L_a h+\nabla\pi=0, \qquad \nabla \cdot h=0 \quad\text{in }Q_R.
\end{align}
Then
\begin{equation}\label{homY}
\Y_a(h;Q_r) \le C(R-r)^{-M}\|h\|_{H^1(Q_R)},
\end{equation}
where \(C,M>0\) are universal constants independent of \(a\).
\end{lemma}

\begin{proof}
For any \(0<r<R\le1\), set $\rho_j=r+\frac{j(R - r)}{8}$ with $\rho_0 = r$ and $\rho_8 = R$.


{\underline{\bf Step 1: removal of the harmonic pressure.}}
Taking the divergence of the equation gives
\begin{align}\label{ine:harmonic pi}
\Delta\pi=0
\quad\text{in }\mathcal D'(Q_R).
\end{align}
Hence \(\pi\) has a smooth harmonic representative. Obviously, 
\begin{align*}
\partial_i\pi=\Delta h_i-a\partial_s h_i,
\qquad i=1,\ldots,5.
\end{align*}
Choose a smooth radial function \(\zeta\in C_c^\infty(B_1)\) satisfying
\begin{align*}
\int_{\mathbb R^5}\zeta(z)\,\mathrm dz=1,
\qquad
\zeta_d(z)=d^{-5}\zeta(z/d).
\end{align*}
For \(x\in Q_{\rho_7}\), take \(d=\frac{R -r}{16}\), so that \(B_d(x)\subset Q_R\).
Since every derivative of \(\pi\) is harmonic, the mean value property implies
\begin{align*}
D^\alpha\partial_i\pi(x) &=-\int_{Q_R} \nabla h_i(z)\cdot \nabla_zD_x^\alpha\zeta_d(x-z)\,\mathrm dz\\* &\quad +a\int_{Q_R} h_i(z)\, \partial_{z_1}D_x^\alpha\zeta_d(x-z)\,\mathrm dz.
\end{align*}
Here \(z_1\) is the longitudinal coordinate. The scaling of the averaging kernel gives
\begin{align*}
\|\nabla D^\alpha\zeta_d\|_{L^2(\mathbb R^5)} \le C_\alpha d^{-|\alpha|-\frac72},
\end{align*}
which means 
\ben\label{ine:partialpi} \nonumber
|D^\alpha\partial_i\pi(x)| &\le& C_\alpha d^{-|\alpha|-\frac72} \left( \|\nabla h\|_{L^2(Q_R)} +a\|h\|_{L^2(Q_R)} \right)\\
&\le& C_\alpha a(R - r)^{-|\alpha|-\frac72} \|h\|_{H^1(Q_R)}.
\een
Let \(\psi\) solve the transverse Dirichlet problem
\begin{align*}
\Delta_y\psi(y)=a^{-1}\partial_s\pi(0,y)
\quad\text{in }B_{\rho_7}^4,
\qquad
\psi=0
\quad\text{on }\partial B_{\rho_7}^4.
\end{align*}
The interior Poisson estimates and $\eqref{ine:partialpi}$ yield
\begin{align*}
\|\psi\|_{L^\infty(B_{\rho_7}^4)} \le C\rho_7^2a^{-1} \|\partial_s\pi(0,\cdot)\|_{L^\infty(B_{\rho_7}^4)} \le C(R - r)^{-M}\|h\|_{H^1(Q_R)}.
\end{align*}

Now for any $(s,y) \in Q_{\rho_6}$, define
\begin{align*}
\Phi(s,y) =-a^{-1}\int_0^s \bigl(\pi(t,y)-\pi(0,0)\bigr)\,\mathrm dt+\psi(y),
\qquad
H=\nabla\Phi.
\end{align*}
By \eqref{ine:harmonic pi}, we have
\begin{align*}
\Delta\Phi(s,y) &=-a^{-1}\partial_s\pi(s,y) -a^{-1}\int_0^s\Delta_y\pi(t,y)\,\mathrm dt +\Delta_y\psi(y)\\
&=-a^{-1}\partial_s\pi(s,y) +a^{-1}\int_0^s\partial_t^2\pi(t,y)\,\mathrm dt +a^{-1}\partial_s\pi(0,y)\\
&=-a^{-1}\partial_s\pi(s,y) +a^{-1}\bigl( \partial_s\pi(s,y)-\partial_s\pi(0,y) \bigr) +a^{-1}\partial_s\pi(0,y)=0.
\end{align*}
Also,
\(a\partial_sH =\nabla(a\partial_s\Phi) =-\nabla\pi.\)
Therefore,
\begin{align}\label{eq:H}
\nabla \cdot  H=0,
\qquad
\Delta H=0,
\qquad
L_a  H+\nabla\pi=0
\quad\text{in }Q_{\rho_6}.
\end{align}
The preceding pressure and Poisson estimates, as detailed in Lemma~\ref{aux:H}, imply 
\ben\label{ine:estimateH}
\|H\|_{H^2(Q_{\rho_6})} \leq C (R - r)^{-M} \|h\|_{H^1(Q_R)},
\een
and the one-dimensional Sobolev inequality on $(-\rho_6,\rho_6)$ gives
\begin{align}\label{ine:estH}
\|\nabla H\|_{L^\infty((-\rho_6,\rho_6);L^2(B_{\rho_6}^4))}^2
&\le C\rho_6^{-1}\|\nabla H\|_{L^2(Q_{\rho_6})}^2+C\|\nabla H\|_{L^2(Q_{\rho_6})}\|\partial_s\nabla H\|_{L^2(Q_{\rho_6})}\notag\\
&\le C(R-r)^{-M}\|h\|_{H^1(Q_R)}^2.
\end{align}

Set \(w=h-H\). Then \eqref{eq:h}, \eqref{eq:H} and \eqref{ine:estimateH} yield
\begin{align}\label{eq:w}
L_a  w=0 \quad\text{in }Q_{\rho_6},
\qquad
\|w\|_{H^1(Q_{\rho_6})} \le C{(R-r)}^{-M}\|h\|_{H^1(Q_R)}.
\end{align}

{\underline{\bf Step 2: full-space correction.}}
Choose
\begin{align*}
\chi\in C_c^\infty(B_{\rho_5}^4),
\qquad
0\le\chi\le1,
\qquad
\chi=1\quad\text{on }B_{\rho_4}^4,
\end{align*}
with
\begin{align*}
\|\nabla_y^k\chi\|_{L^\infty(\mathbb R^4)} \le C{(R-r)}^{-k}.
\end{align*}
Put \(I=(-\rho_5,\rho_5)\). Since \(\chi\) is independent of \(s\), the product rule and \eqref{eq:w} give
\begin{align*}
L_a (\chi w) =-2\nabla_y\chi\cdot\nabla_yw -(\Delta_y\chi)w =:F \quad\text{in } I\times\mathbb R^4.
\end{align*}
Extend \(F\) by zero outside \(I\times B_{\rho_5}^4\). Then \eqref{eq:w} yields
\begin{align*}
\|F\|_{L^2(\mathbb R^5)} &\le C\left( {(R-r)}^{-1}\|\nabla_yw\|_{L^2(Q_{\rho_5})} +{(R-r)}^{-2}\|w\|_{L^2(Q_{\rho_5})} \right)\\
&\le C{(R-r)}^{-M}\|h\|_{H^1(Q_R)}.
\end{align*}

Letting \(Z=L_a ^{-1}F\), its Fourier representation is
\begin{align*}
\widehat Z(\xi)=D_a^{-1}\widehat F(\xi),
\qquad
D_a=|\xi|^2+ia\tau,
\qquad
\xi=(\tau,\xi').
\end{align*}
A direct calculation shows that
\begin{align}\label{ine:bddD}
|D_a|^{-1}\le|\xi|^{-2},
\qquad
|\xi|\,|D_a|^{-1}\le|\xi|^{-1},
\qquad
|\tau|\,|\xi|\,|D_a|^{-1}\le1.
\end{align}
By Plancherel's theorem, \eqref{ine:bddD} and Lemma \ref{lem:HLSine}, there holds
\begin{align*}
&\|Z\|_{L^2(\mathbb R^5)} \le\||\xi|^{-2} \widehat F\|_{L^2(\mathbb R^5)} = \|(-\Delta)^{-1} F\|_{L^2(\mathbb R^5)} \le C\|F\|_{L^{\frac{10}{9}}(\mathbb R^5)}\le C\|F\|_{L^2(\mathbb R^5)},\\
&\|\nabla Z\|_{L^2(\mathbb R^5)} \leq \||\xi|^{-1} \widehat F\|_{L^2(\mathbb R^5)} \le\|(-\Delta)^{-\frac12}F\|_{L^2(\mathbb R^5)} \le C\|F\|_{L^{\frac{10}{7}}(\mathbb R^5)}\le C\|F\|_{L^2(\mathbb R^5)},\\
&\|\partial_s\nabla Z\|_{L^2(\mathbb R^5)}^2 \le\int_{\mathbb R^5}\frac{\tau^2|\xi|^2}{|D_a|^2} |\widehat F(\xi)|^2\,\mathrm d\xi \le\|F\|_{L^2(\mathbb R^5)}^2.
\end{align*}
Thus
\ben\label{ine:estZh}
\|Z\|_{H^1(\mathbb R^5)} +\|\partial_s\nabla Z\|_{L^2(\mathbb R^5)} \le C{(R-r)}^{-M}\|h\|_{H^1(Q_R)}.
\een
The one-dimensional Sobolev inequality and \eqref{ine:estZh} give
\begin{align}\label{ine:estZ}
\|\nabla Z\|_{L^\infty(\mathbb R;L^2(\mathbb R^4))}^2 \le 2\|\nabla Z\|_{L^2(\mathbb R^5)} \|\partial_s\nabla Z\|_{L^2(\mathbb R^5)}\le C{(R-r)}^{-M}\|h\|_{H^1(Q_R)}^2.
\end{align}

Now we define 
\begin{align*}
W=\chi w-Z
\quad\text{on }I\times\mathbb R^4.
\end{align*}
Since \(L_a  Z=F\) in the full space,
\begin{align}\label{eq:W}
L_a  W=0
\quad\text{in }I\times\mathbb R^4.
\end{align}
Moreover, \eqref{eq:w} and \eqref{ine:estZh} yield
\begin{align}\label{ine:estW}
\|W\|_{H^1(I\times\mathbb R^4)} &\le \|\chi w\|_{H^1(I\times\mathbb R^4)} +\|Z\|_{H^1(\mathbb R^5)}\le C{(R-r)}^{-M}\|h\|_{H^1(Q_R)}.
\end{align}

{\underline{\bf Step 3: convexity of the slice gradient energy.}}
Let
\begin{align*}
\widehat W(s,\xi')=\mathcal F_yW(s,\xi'), \qquad k=|\xi'|.
\end{align*}
Taking the Fourier transform on both sides of equation \eqref{eq:W}, we obtain
\begin{align*}
-\partial_s^2\widehat W +a\partial_s\widehat W+k^2\widehat W=0.
\end{align*}
For almost every \(\xi'\), its solution has the form
\begin{align*}
\widehat W(s,\xi') =c_+(\xi')e^{ \lambda_+s} +c_-(\xi')e^{ \lambda_-s},
\qquad
 \lambda_\pm=\frac{a\pm\sqrt{a^2+4k^2}}2.
\end{align*}
The characteristic roots satisfy
\begin{align*}
 \lambda_++ \lambda_-=a,
\qquad
\lambda_+ \lambda_-=-k^2.
\end{align*}
Consequently, the complete gradient energy at each frequency satisfies the exact identity
\begin{align*}
|\partial_s\widehat W|^2+k^2|\widehat W|^2 ={}&( \lambda_+^2+k^2)|c_+|^2e^{2 \lambda_+s} +(\lambda_-^2+k^2)|c_-|^2e^{2 \lambda_-s}\\
&+2( \lambda_+ \lambda_-+k^2) \operatorname{Re}(c_+\cdot\overline{c_-}) e^{( \lambda_++ \lambda_-)s}\\
={}&( \lambda_+^2+k^2)|c_+|^2e^{2 \lambda_+s} +( \lambda_-^2+k^2)|c_-|^2e^{2 \lambda_-s}.
\end{align*}
Clearly, \(|\partial_s\widehat W|^2+k^2|\widehat W|^2\) is a convex function of \(s\).
Now we define
\begin{align*}
q_W(s)=\int_{\mathbb R^4} \left(|\partial_s\widehat W(s,\xi')|^2 +|\xi'|^2|\widehat W(s,\xi')|^2\right)\,\mathrm d\xi'.
\end{align*}
For \([s-d,s+d]\subset I\), by Jensen's inequality, Tonelli's theorem, Plancherel's theorem and \eqref{ine:estW}, we have 
\begin{align*}
q_W(s) &\le \frac1{2d}\int_{s-d}^{s+d}q_W(t)\,\mathrm dt = \frac1{2d}\int_{s-d}^{s+d}\int_{\mathbb R^4} |\nabla W(t,y)|^2\,\mathrm dy\,\mathrm dt \\
&\le \frac{\|\nabla W\|_{L^2(I\times\mathbb R^4)}^2}{2d}\le \frac C{2d} {(R-r)}^{-M}\|h\|_{H^1(Q_R)}^2.
\end{align*}
Taking \(d=\frac1{16}(R-r)\) for \(s\in(-\rho_4,\rho_4)\), we obtain
\begin{align}
\sup_{|s|<\rho_4}q_W(s) &\le C{(R-r)}^{-M}\|h\|_{H^1(Q_R)}^2, \\
\|\nabla W\|_{L^\infty((-\rho_4,\rho_4);L^2(\mathbb R^4))} &\le C{(R-r)}^{-M}\|h\|_{H^1(Q_R)}. \label{ine:nablaW}
\end{align}

On $Q_{\rho_4}$, $\chi=1$, so $h=H+W+Z$. Combining \eqref{ine:nablaW}, \eqref{ine:estZ} and \eqref{ine:estH}, we have 
\begin{align*}
\|\nabla h\|_{L^\infty((-\rho_4,\rho_4); L^2(B_{\rho_4}^4))} &\le \|\nabla H\|_{L^\infty((-\rho_4,\rho_4); L^2(B_{\rho_4}^4))}\\
&+ \|\nabla W\|_{L^\infty((-\rho_4,\rho_4); L^2(\mathbb R^4))}+ \|\nabla Z\|_{L^\infty(\mathbb R; L^2(\mathbb R^4))} \\
&\le C{(R-r)}^{-M}\|h\|_{H^1(Q_R)}.
\end{align*}
In particular,
\(\|\nabla h\|_{L^\infty((-r,r);L^2(B_r^4))} \le C{(R-r)}^{-M}\|h\|_{H^1(Q_R)}.\) By the Sobolev inequality, there holds
\begin{align}\label{ine:Yh1}
\|h\|_{L^\infty((-\rho_4,\rho_4); L^2(B_{\rho_4}^4))} &\le C\left( \rho_4^{-\frac12}\|h\|_{L^2(Q_{\rho_4})} +\rho_4^{\frac12}\|\partial_sh\|_{L^2(Q_{\rho_4})} \right)\le C{(R-r)}^{-\frac12}\|h\|_{H^1(Q_R)},
\end{align}
which means that 
\begin{align}\label{ine:Yh2} \nonumber
\|h\|_{L^\infty((-r,r);L^4(B_r^4))} &\le \|h\|_{L^\infty((-\rho_4,\rho_4); L^4(B_{\rho_4}^4))}\\ \nonumber
&\le C \|\nabla_yh\|_{L^\infty((-\rho_4,\rho_4); L^2(B_{\rho_4}^4))}+ C\rho_4^{-1} \|h\|_{L^\infty((-\rho_4,\rho_4); L^2(B_{\rho_4}^4))}\\
&\le C{(R-r)}^{-M}\|h\|_{H^1(Q_R)}.
\end{align}

{\underline{\bf Step 4: the second-derivative estimate.}}
Since \(L_a \) has constant coefficients, each scalar function
\(g=\partial_iw_j\) satisfies
\begin{align*}
L_a  g=0
\quad\text{in }Q_{\rho_6}.
\end{align*}
Choose
\begin{align*}
\eta\in C_c^\infty(Q_{\rho_4}),
\qquad
0\le\eta\le1,
\qquad
\eta=1\quad\text{on }Q_r,
\end{align*}
such that
\(\|\nabla\eta\|_{L^\infty(\mathbb R^5)} +\|\partial_s\eta\|_{L^\infty(\mathbb R^5)} \le C{(R-r)}^{-1}.\)
For fixed $a$, interior elliptic regularity gives $g\in C^\infty(Q_{\rho_6})$. Testing its equation with $g\eta^2$ and integrating by parts gives
\begin{align*}
\int_{Q_{\rho_4}}|\nabla g|^2\eta^2\,\mathrm dx
&=-2\int_{Q_{\rho_4}}
g\eta\nabla g\cdot\nabla\eta\,\mathrm dx
+a\int_{Q_{\rho_4}}|g|^2\eta\partial_s\eta\,\mathrm dx.
\end{align*}
Young's inequality therefore yields
\begin{align*}
\int_{Q_{\rho_4}}|\nabla g|^2\eta^2\,\mathrm dx &\le\frac12 \int_{Q_{\rho_4}}|\nabla g|^2\eta^2\,\mathrm dx+ C({(R-r)}^{-2}+a{(R-r)}^{-1}) \int_{Q_{\rho_4}}|g|^2\,\mathrm dx.
\end{align*}
Since \(a\ge1\) and \(0<{(R-r)}\le1\),
\begin{align}\label{ine:D2h} \nonumber
a^{-\frac12}\|\nabla^2w\|_{L^2(Q_r)} &\le C\left( a^{-\frac12}{(R-r)}^{-1}+{(R-r)}^{-\frac12} \right) \|\nabla w\|_{L^2(Q_{\rho_4})} \\
&\le C{(R-r)}^{-1}\|w\|_{H^1(Q_{\rho_6})} \le C{(R-r)}^{-M}\|h\|_{H^1(Q_R)}.
\end{align}
Since $h=w+H$ in $Q_r$, by \eqref{ine:estimateH} and \eqref{ine:D2h}, there holds
\begin{align}\label{ine:Yh3}
a^{-\frac12}\|\nabla^2h\|_{L^2(Q_r)} &\le a^{-\frac12}\|\nabla^2w\|_{L^2(Q_r)} +a^{-\frac12}\|\nabla^2H\|_{L^2(Q_r)}\le C{(R-r)}^{-M}\|h\|_{H^1(Q_R)}.
\end{align}

Combining \eqref{ine:nablaW} with \eqref{ine:Yh1}, \eqref{ine:Yh2} and \eqref{ine:Yh3}, we obtain
\begin{align*}
\begin{aligned}
\Y_a(h;Q_r) ={}&\|h\|_{L^\infty((-r,r);L^4(B_r^4))} +\|\nabla h\|_{L^\infty((-r,r);L^2(B_r^4))}+a^{-\frac12}\|\nabla^2h\|_{L^2(Q_r)}\\ \le{}&C{(R-r)}^{-M}\|h\|_{H^1(Q_R)}. 
\end{aligned}
\end{align*}
Here \(M\) depends on the order of differentiation, which is fixed
in the present proof. This proves \eqref{homY}.
\end{proof}

\begin{lemma}\label{linearized}
Let \(\frac12\le r<R\le1\), \(a\ge1\), and let \(U\) be smooth and divergence free in \(Q_R\).
Suppose that
\begin{align*}
h\in H^1(Q_R),
\qquad
\pi\in L^{\frac53}_{\mathrm{loc}}(Q_R),
\end{align*}
and
\begin{align*}
L_a h+ \nabla \cdot(U\otimes h+h\otimes U)+\nabla\pi=0,
\qquad
\nabla \cdot h=0
\quad\text{in }\mathcal D'(Q_R).
\end{align*}
There exists a universal constant \(\delta_{\mathrm{lin}}>0\)
such that, if \(\Y_a(U;Q_R)\le\delta_{\mathrm{lin}}\), then
\begin{equation}\label{linY}
\Y_a(h;Q_r) \le C(R-r)^{-M}\|h\|_{H^1(Q_R)},
\end{equation}
where \(C,M>0\) are universal constants independent of \(a\) and of higher smooth norms of \(U\).
\end{lemma}

\begin{proof}
Without loss of generality, we require only \(\delta_{\mathrm{lin}}\le1\).

{\underline{\bf Step 1: qualitative interior regularity.}}
For fixed \(a\), write
\begin{align*}
-\Delta h+\nabla\pi =-a\partial_sh-U\cdot\nabla h-h\cdot\nabla U.
\end{align*}
Since \(U\) is smooth, standard interior Stokes regularity and bootstrapping yield
\begin{align}\label{ine:Ginfty}
h,\pi\in C^\infty(Q_R), \qquad \Y_a(h;Q_\rho)<\infty \quad\text{for every }0<\rho<R.
\end{align}

{\underline{\bf Step 2: reduction to a homogeneous equation.}}
Fix
\begin{align*}
r\le t<\rho<R, \qquad \sigma=\frac{t+\rho}{2}.
\end{align*}
Choose a cutoff satisfying
\begin{align*}
\eta\in C_c^\infty(Q_\rho), \qquad 0\le\eta\le1, \qquad \eta=1\quad\text{on }Q_\sigma, \qquad
\|\nabla\eta\|_{L^\infty(\mathbb R^5)} \le C(\rho-t)^{-1}.
\end{align*}
Define the tensor \(T=\eta(U\otimes h+h\otimes U)\) on \(Q_\rho\), extend it by zero to \(\mathbb R^5\), and set \(Z=S_a  T=\mathcal R_a \nabla \cdot  T.\)
By H\"older's inequality, there holds
\begin{align*}
\|T\|_{L^2(\mathbb R^5)}^2 &\le 4\int_{-\rho}^{\rho} \|U(s,\cdot)\|_{L^4(B_\rho^4)}^2 \|h(s,\cdot)\|_{L^4(B_\rho^4)}^2 \,\mathrm ds \le 4\delta_{\mathrm{lin}}^2 \|h\|_{L^2((-\rho,\rho);L^4(B_\rho^4))}^2.
\end{align*}
For $\tfrac12\le R\le1$, the Sobolev inequality gives
\begin{equation*}
\begin{aligned}
\|h\|_{L^2((-R,R);L^4(B_R^4))}^2 &=\int_{-R}^{R}\|h(s,\cdot)\|_{L^4(B_R^4)}^2\,\mathrm ds\\
&\le C\int_{-R}^{R}\int_{B_R^4} \left(|\nabla_y h(s,y)|^2+|h(s,y)|^2\right)\,\mathrm dy\,\mathrm ds
\le C\|h\|_{H^1(Q_R)}^2.
\end{aligned}
\end{equation*}
Then we have
\ben\label{ine:estT}
\|T\|_{L^2(\mathbb R^5)}^2 \leq C \delta_{\mathrm{lin}}^2 \|h\|_{H^1(Q_R)}^2.
\een
Since \(T\) is supported in \(Q_1\), \eqref{compactT} and \eqref{ine:estT} give
\begin{align}\label{ine:estimateZ}
\|Z\|_{H^1(\mathbb R^5)} \le C\|T\|_{L^2(\mathbb R^5)} \le C\delta_{\mathrm{lin}}\|h\|_{H^1(Q_R)}.
\end{align}

Since the Leray projection makes \(Z\) divergence free, we can choose a distribution \(p_Z\) such that \(\nabla p_Z=(I-\mathbb P)\nabla \cdot  T.\)
Then
\begin{align*}
 L_a  Z+\nabla p_Z=\nabla \cdot  T, \qquad \nabla \cdot  Z=0 \quad\text{in }\mathbb R^5.
\end{align*}
Because \(\eta=1\) on \(Q_\sigma\), the equations for \(h\) and \(Z\) imply
\begin{align*}
 L_a (h+Z)+\nabla(\pi+p_Z)=0, \qquad \nabla \cdot (h+Z)=0 \quad\text{in }Q_\sigma.
\end{align*}
Applying Lemma~\ref{homogeneous} on \(Q_t\subset Q_\sigma\) and \eqref{ine:estimateZ}, we obtain
\begin{align}\label{ine:estYZ+h} \nonumber
\Y_a(h+Z;Q_t) &\le C(\sigma-t)^{-M}\|h+Z\|_{H^1(Q_\sigma)} \le C(\rho-t)^{-M} \left( \|h\|_{H^1(Q_\rho)} +\|Z\|_{H^1(\mathbb R^5)} \right)\\
&\le C(\rho-t)^{-M}\|h\|_{H^1(Q_R)}.
\end{align}

{\underline{\bf Step 3: estimate of the correction in \(\Y_a\).}}
Set \(F=\nabla \cdot  T\). Since \(U\) and \(h\) are divergence free, the product rule gives
\begin{align*}
F =\eta\bigl((U\cdot\nabla)h+(h\cdot\nabla)U\bigr) +(U\otimes h+h\otimes U)\nabla\eta.
\end{align*}
We estimate $F$ term by term. By H\"older's inequality and the Sobolev inequality, we obtain
\begin{align}\label{ine:estimateF1} \nonumber
\|(U\cdot\nabla)h\|_{L^2(Q_\rho)} &\le \|U\|_{L^\infty((-\rho,\rho);L^4(B_\rho^4))} \|\nabla h\|_{L^2((-\rho,\rho);L^4(B_\rho^4))}\\ \nonumber
&\le C\delta_{\mathrm{lin}} \left( \|\nabla^2h\|_{L^2(Q_\rho)} +\|\nabla h\|_{L^2(Q_\rho)} \right) \\
&\le C\sqrt a\,\delta_{\mathrm{lin}}\Y_a(h;Q_\rho) +C\delta_{\mathrm{lin}}\|h\|_{H^1(Q_R)}.
\end{align}
For the second term, the same method gives
\begin{align*}
\|(h\cdot\nabla)U\|_{L^2(Q_\rho)} &\le \|h\|_{L^\infty((-\rho,\rho);L^4(B_\rho^4))} \|\nabla U\|_{L^2((-\rho,\rho);L^4(B_\rho^4))}\\
&\le C\Y_a(h;Q_\rho) \left( \|\nabla^2U\|_{L^2(Q_\rho)} +\|\nabla U\|_{L^2(Q_\rho)} \right).
\end{align*}
Since $\Y_a(U;Q_R) \leq \delta_{\mathrm{lin}}$ and
\begin{align*}
\|\nabla U\|_{L^2(Q_\rho)}^2 &=\int_{-\rho}^{\rho} \|\nabla U(s,\cdot)\|_{L^2(B_\rho^4)}^2 \,\mathrm ds\le 2\rho\, \|\nabla U\|_{L^\infty((-\rho,\rho);L^2(B_\rho^4))}^2,
\end{align*}
there holds
\ben\label{ine:estimateF2}
\|(h\cdot\nabla)U\|_{L^2(Q_\rho)} \le C\sqrt a\,\delta_{\mathrm{lin}}\Y_a(h;Q_\rho).
\een
For the last term, noting that $|\nabla \eta| \leq C (\rho - t)^{-1}$, there holds
\begin{align}\label{ine:estimateF3} \nonumber
\|(U\otimes h+h\otimes U)\nabla\eta\|_{L^2(\mathbb R^5)} &\le C(\rho-t)^{-1} \|U\|_{L^\infty((-\rho,\rho);L^4(B_\rho^4))} \|h\|_{L^2((-\rho,\rho);L^4(B_\rho^4))}\\
&\le C(\rho-t)^{-1}\delta_{\mathrm{lin}} \|h\|_{H^1(Q_R)}.
\end{align}
Combining \eqref{ine:estimateF1}, \eqref{ine:estimateF2}, \eqref{ine:estimateF3} and using \(\rho-t\le1\), we have
\begin{align*}
\|F\|_{L^2(\mathbb R^5)} \le C\sqrt a\,\delta_{\mathrm{lin}}\Y_a(h;Q_\rho) +C(\rho-t)^{-1}\delta_{\mathrm{lin}}\|h\|_{H^1(Q_R)}.
\end{align*}
Since \(Z=\mathcal R_a  F\), \eqref{RF} implies
\begin{align}\label{ine:estYZ}\nonumber
\Y_a(Z)&\le Ca^{-\frac12}\|F\|_{L^2(\mathbb R^5)} \le C\delta_{\mathrm{lin}}\Y_a(h;Q_\rho) +Ca^{-\frac12}(\rho-t)^{-1}\delta_{\mathrm{lin}}\|h\|_{H^1(Q_R)}\\
&\le C\delta_{\mathrm{lin}}\Y_a(h;Q_\rho) +C(\rho-t)^{-1}\delta_{\mathrm{lin}}\|h\|_{H^1(Q_R)}.
\end{align}
Noting that $h=(h+Z)-Z$, we combine \eqref{ine:estYZ} and \eqref{ine:estYZ+h}, increasing the fixed exponent $M$ if necessary, to obtain
\begin{align}\label{Yhole} \nonumber
\Y_a(h;Q_t) &\le\Y_a(h+Z;Q_t)+\Y_a(Z;Q_t)\\
&\le C\delta_{\mathrm{lin}}\Y_a(h;Q_\rho) +C(\rho-t)^{-M}\|h\|_{H^1(Q_R)}.
\end{align}

{\underline{\bf Step 4: iterative discussion.}}
Set 
\begin{align*}
R_\ast = \frac {r+R}2, \qquad  G(t)=\Y_a(h;Q_t).
\end{align*}
By \eqref{ine:Ginfty} and \eqref{Yhole}, we have
\begin{align*}
G(t)\le C\delta_{\mathrm{lin}}G(\rho) +C\|h\|_{H^1(Q_R)}(\rho-t)^{-M}, \qquad r\le t<\rho\le R_\ast.
\end{align*}
Choose \(0<\delta_{\mathrm{lin}}\le1\) so that \(C\delta_{\mathrm{lin}}2^M\le\frac12\). Lemma~\ref{lem:Gia} yields
\begin{align*}
\Y_a(h;Q_r) \le C(R_\ast-r)^{-M}\|h\|_{H^1(Q_R)} \le C(R-r)^{-M}\|h\|_{H^1(Q_R)}.
\end{align*}
This proves the lemma.
\end{proof}

\section{Quadratic improvement of the perturbation} 

\subsection{Smooth corrections and the relative energy inequality}

\par We first incorporate a small linearized perturbation into an
exact smooth Navier--Stokes background by a quadratic correction.
We then derive a relative energy inequality with a smooth
comparison field.

We extend the divergence-free fields from a smaller cylinder using Lemma \ref{extension} in the appendix.

\begin{lemma}\label{strongcorrect}
Let \(a\ge1\), and suppose that \(ae_1+U\) is a smooth divergence-free solution of the stationary Navier-Stokes equations in \(Q_{\frac12}\).
Let \(h\) be a smooth divergence-free solution of its homogeneous linearization there, so that
\begin{align*}
\begin{aligned}
L_a U+(U\cdot\nabla)U+\nabla p_U&=0,\\
L_a h+(U\cdot\nabla)h+(h\cdot\nabla)U+\nabla\pi&=0
\end{aligned}
\qquad\text{in }Q_{\frac12}.
\end{align*}
Suppose
\begin{align*}
\Y_a(U;Q_{\frac12})\le\delta, \qquad \Y_a(h;Q_{\frac12})\le C_0\varepsilon,
\end{align*}
where \(C_0>0\) is fixed. If \(\delta+C_0\varepsilon\) is sufficiently small, there exists a divergence-free field \(k\) such that \(\Y_a(k)\le C C_0^2\varepsilon^2\) and \(ae_1+U+h+k\) is a smooth exact Navier-Stokes solution in \(Q_{\frac14}\).
The smallness threshold and \(C\) are independent of \(a\).
\end{lemma}

\begin{proof}

{\underline{\bf Step 1: extension and the fixed-point equation.}}
By Lemma \ref{extension}, there are divergence-free extensions $\widehat U$ and $\widehat h$ such that
\begin{align*}
\widehat U=U,\qquad \widehat h=h \quad\text{in }Q_{\frac14},
\end{align*}
and
\begin{align*}
\Y_a(\widehat U)\le C_{\mathrm{ext}}\delta, \qquad \Y_a(\widehat h)\le C_{\mathrm{ext}}C_0\varepsilon.
\end{align*}
The constant \(C_{\mathrm{ext}}\) is fixed.

Define
\begin{align*}
\begin{aligned}
\mathcal T(k)=-\mathcal R_a \bigl[ (\widehat U+\widehat h)\cdot\nabla k +k\cdot\nabla(\widehat U+\widehat h)+k\cdot\nabla k+\widehat h\cdot\nabla\widehat h
\bigr].
\end{aligned}
\end{align*}
We seek a fixed point of \(\mathcal T\) in the Banach space of
divergence-free fields \(k\in\dot H^2(\mathbb R^5)\) with finite
\(\Y_a\)-norm, equipped with this norm. By H\"older's inequality and the Sobolev inequality, there holds
\begin{align*}
\begin{aligned}
\|V\cdot\nabla W\|_{L^2(\mathbb R^5)} &\le \|V\|_{L^\infty(\mathbb R;L^4(\mathbb R^4))} \|\nabla W\|_{L^2(\mathbb R;L^4(\mathbb R^4))}\\
&\le C \|V\|_{L^\infty(\mathbb R;L^4(\mathbb R^4))} \|\nabla^2W\|_{L^2(\mathbb R^5)} \le C\sqrt a\,\Y_a(V)\Y_a(W).
\end{aligned}
\end{align*}

%

{\underline{\bf Step 2: contraction.}} Fix a sufficiently large universal constant \(B\ge1\) and assume
\begin{align*}
\delta+C_0\varepsilon\le(8BC_{\mathrm{ext}})^{-1}.
\end{align*}
Then
\begin{align}\label{ine:smallYUh}
2B\bigl(\Y_a(\widehat U)+\Y_a(\widehat h)\bigr) +4B^2\Y_a(\widehat h)^2 < \frac12.
\end{align} 
On the ball $\Y_a(k)\le2B\Y_a(\widehat h)^2$, Corollary \ref{cor:bilinear-estimate} and \eqref{ine:smallYUh} give 
\ben\label{ine:YakBh}\nonumber
\Y_a(\mathcal T(k)) &\leq& 2B\bigl(\Y_a(\widehat U)+\Y_a(\widehat h)\bigr)\Y_a(k)+B\Y_a(k)^2+B\Y_a(\widehat h)^2\\
&\leq& \frac12\Y_a(k)+B\Y_a(\widehat h)^2. 
\een
On the closed ball \(\Y_a(k)\le2B\Y_a(\widehat h)^2\), \eqref{ine:YakBh} yields
\begin{align*}
\Y_a(\mathcal T(k)) \leq 2B\Y_a(\widehat h)^2.
\end{align*}
Similarly, 
\begin{align*}
\Y_a\bigl(\mathcal T(k_1)-\mathcal T(k_2)\bigr) &\le B\bigl[ 2\Y_a(\widehat U)+2\Y_a(\widehat h) +\Y_a(k_1)+\Y_a(k_2)\bigr]\Y_a(k_1-k_2) \\
&\le\frac12\Y_a(k_1-k_2).
\end{align*}
The Banach fixed-point theorem therefore gives \(k=\mathcal T(k)\), with
\begin{align*}
\Y_a(k)\le2B\Y_a(\widehat h)^2 \le2BC_{\mathrm{ext}}^2C_0^2\varepsilon^2.
\end{align*}

{\underline{\bf Step 3: recovery of the exact equation.}}
The fixed-point identity implies
\begin{align*}
\begin{aligned}
L_a k=-\mathbb P \bigl[ (\widehat U+\widehat h)\cdot\nabla k +k\cdot\nabla(\widehat U+\widehat h)+k\cdot\nabla k+\widehat h\cdot\nabla\widehat h \bigr].
\end{aligned}
\end{align*}
So a pressure \(\pi_k\) can be chosen with
\begin{align*}
\begin{aligned}
\nabla\pi_k=-(I-\mathbb P )\bigl[ (\widehat U+\widehat h)\cdot\nabla k +k\cdot\nabla(\widehat U+\widehat h)+k\cdot\nabla k+\widehat h\cdot\nabla\widehat h
\bigr].
\end{aligned}
\end{align*}
Since the extensions equal the original fields in \(Q_{\frac14}\), we obtain that on $Q_\frac14$, 
\begin{align*}
\begin{aligned}
&L_a U+(U\cdot\nabla)U+\nabla p_U=0,\\
&L_a h+(U\cdot\nabla)h+(h\cdot\nabla)U+\nabla\pi=0,\\
&L_a k+((U+h)\cdot\nabla)k +(k\cdot\nabla)(U+h) +(k\cdot\nabla)k+\nabla\pi_k =-(h\cdot\nabla)h.
\end{aligned}
\end{align*}
Adding these equations and collecting all quadratic terms
gives
\begin{align*}
\begin{aligned}
L_a(U+h+k) +((U+h+k)\cdot\nabla)(U+h+k) +\nabla(p_U+\pi+\pi_k)=0,
\end{aligned}
\end{align*}
Thus \(v=ae_1+U+h+k\) solves the stationary Navier--Stokes
equations in \(Q_{\frac14}\). Since \(v\in H^2_{\rm loc}(Q_{\frac14})
\subset L^{10}_{\rm loc}(Q_{\frac14})\), interior Stokes bootstrapping
gives \(v\in C^\infty(Q_{\frac14})\).

\end{proof}

\begin{lemma}\label{relativeLemma}
Let \((u,p)\) be a suitable weak solution in \(Q_{\frac78}\), and let \(V,P,F\) be smooth functions with
\begin{align}\label{eq:VPF}
\nabla \cdot  V=0,\qquad -\Delta V+(V\cdot\nabla)V+\nabla P=F.
\end{align}
Set \(z=u-V\) and \(q=p-P\). For any nonnegative \(\phi\in C_c^\infty(Q_\frac78)\), there holds
\begin{align}
\int_{{Q_\frac78}}|\nabla z|^2\phi \, dx  \le{}&\frac12\int_{{Q_\frac78}}|z|^2\Delta\phi \, dx  +\int_{{Q_\frac78}}\left(\frac12|z|^2u+qz\right)\cdot\nabla\phi \, dx \notag\\
&-\int_{{Q_\frac78}}(z\cdot\nabla V)\cdot z\,\phi \, dx  -\int_{{Q_\frac78}} F\cdot z\,\phi \, dx .\label{generalRelative}
\end{align}
\end{lemma}
\begin{proof}
Testing \eqref{NS} with $V\phi$ gives
\begin{align*}
\int_{{Q_\frac78}}\nabla u:\nabla V\,\phi \, dx &= \int_{{Q_\frac78}}(u\cdot\nabla V)\cdot u\,\phi \, dx  +\int_{{Q_\frac78}}(u\cdot V)u\cdot\nabla\phi \, dx  \\
&+\int_{{Q_\frac78}}pV\cdot\nabla\phi \, dx  -\int_{{Q_\frac78}}\partial_ju_i V_i\partial_j\phi \, dx .
\end{align*}
Testing \eqref{eq:VPF} with $u\phi$ gives
\begin{align*}
\int_{{Q_\frac78}}\nabla u:\nabla V\,\phi \, dx &=\int_{{Q_\frac78}}F\cdot u\,\phi \, dx  -\int_{{Q_\frac78}}(V\cdot\nabla V)\cdot u\,\phi \, dx \\
&+\int_{{Q_\frac78}}P u\cdot\nabla\phi \, dx  -\int_{{Q_\frac78}}\partial_jV_i u_i\partial_j\phi \, dx .
\end{align*}
Adding these identities and integrating the last two derivative terms yields
\begin{align}
2\int_{{Q_\frac78}}\nabla u:\nabla V\,\phi \, dx ={}&\int_{{Q_\frac78}}(z\cdot\nabla V)\cdot u\,\phi \, dx  +\int_{{Q_\frac78}}(u\cdot V)u\cdot\nabla\phi \, dx \notag\\
&+\int_{{Q_\frac78}}(pV+Pu)\cdot\nabla\phi \, dx  +\int_{{Q_\frac78}}F\cdot u\,\phi \, dx  +\int_{{Q_\frac78}}(u\cdot V)\Delta\phi \, dx .\label{crossIdentity}
\end{align}
Testing \eqref{eq:VPF} with $V\phi$ gives
\begin{align}\label{eq:local energy V}
\int_{{Q_\frac78}}|\nabla V|^2\phi \, dx  =\frac12\int_{{Q_\frac78}}|V|^2\Delta\phi \, dx  +\int_{{Q_\frac78}}\left(\frac12|V|^2+P\right)V\cdot\nabla\phi \, dx  +\int_{{Q_\frac78}}F\cdot V\,\phi \, dx .
\end{align}
 Using \(|\nabla z|^2=|\nabla u|^2+|\nabla V|^2-2\nabla u:\nabla V\), by \eqref{eq:local energy}, \eqref{eq:local energy V} and \eqref{crossIdentity}, we obtain \eqref{generalRelative}.
\end{proof}

\subsection{The one-step quadratic improvement}

We combine the preceding estimates to obtain a quadratic
improvement of the weak perturbation while preserving an
exact smooth Navier--Stokes background.
The improvement is uniform for \(a\ge1\).

%

\begin{proposition}\label{improvement}
There exist universal constants \(\bar\delta,\bar\varepsilon,C>0\) such that the following statement holds for every \(a\ge1\).
Suppose that \((u,p)\) is a suitable weak solution of the stationary incompressible Navier-Stokes equations in \(Q_1\),
and that \(u=ae_1+U+w,\) where \(ae_1+U\) is a smooth exact solution in \(Q_1\). Assume
\begin{align}\label{ass:keyprop}
\Y_a(U;Q_1)\le\delta\le\bar\delta,
\qquad
\|w\|_{H^1(Q_1)}\le\varepsilon\le\bar\varepsilon.
\end{align}
Then there exists a smooth exact background \(ae_1+U_+\) in \(Q_{\frac14}\) such that
\begin{equation}\label{improve}
\Y_a(U_+;Q_{\frac14})\le\delta+C\varepsilon,
\qquad
\|u-ae_1-U_+\|_{H^1(Q_{\frac14})} \le C\varepsilon^2.
\end{equation}
\end{proposition}

\begin{proof}
Without loss of generality, we assume that \(\delta\le1\) and \(\varepsilon\le1\).

{\underline{\bf Step 1: construction of the weak corrector.}}
By Lemma \ref{extension}, there exists a divergence-free extension \(\widetilde U\), compactly supported in \(Q_1\), such that
\begin{align*}
\widetilde U=U \quad\text{on }{Q_\frac78},
\qquad
\Y_a(\widetilde U)\le C_{\mathrm{ext}}\delta.
\end{align*}
Choose \(\chi\in C_c^\infty(Q_{\frac{15}{16}})\) satisfying $0\le\chi\le1$, and $\chi=1$ on ${Q_\frac78}$. Extend \(\chi w\otimes w\) by zero to \(\mathbb R^5\). 
By Sobolev embedding, we have
\begin{align}\label{ine:w103}
\|\chi w\otimes w\|_{L^{\frac53}(\mathbb R^5)} \le \|w\|_{L^{\frac{10}{3}}(Q_1)}^2 \le C\|w\|_{H^1(Q_1)}^2 \le C\varepsilon^2.
\end{align}
On the other hand, the Sobolev estimate gives
\begin{equation}\label{ine:wL2L4}
\begin{aligned}
\|\chi w\otimes w\|_{L^1(\mathbb R;L^2(\mathbb R^4))} &= \int_{\mathbb R} \left( \int_{\mathbb R^4} |\chi(s,y)w(s,y)\otimes w(s,y)|^2 \,\mathrm dy
\right)^{\frac12} \,\mathrm ds\\
&\le \int_{-1}^{1} \|w(s,\cdot)\|_{L^4(B_1^4)}^2 \,\mathrm ds = \|w\|_{L^2((-1,1);L^4(B_1^4))}^2 \le C\|w\|_{H^1(Q_1)}^2 \\
&\le C\varepsilon^2.
\end{aligned}
\end{equation}
Define
\begin{align*}
z_0=- S_a (\chi w\otimes w),
\qquad
z=(I+B_{\widetilde U})^{-1}z_0,
\end{align*}
where
\(B_{\widetilde U}f =  S_a (\widetilde U\otimes f+f\otimes\widetilde U).\) Thus \eqref{CZ} and \eqref{driftL2} imply
\(\sqrt a\,\|z_0\|_{L^2(\mathbb R^5)} + \|\nabla z_0\|_{L^{\frac53}(\mathbb R^5)} \le C\varepsilon^2.\)
Choose \(\bar\delta\) sufficiently small that Lemma \ref{cross} applies. By the homogeneous Sobolev embedding \(\dot W^{1,\frac53}(\mathbb R^5) \hookrightarrow L^{\frac52}(\mathbb R^5)\), there holds
\begin{equation}\label{weakz}
\sqrt a\,\|z\|_{L^2(\mathbb R^5)} + \|\nabla z\|_{L^{\frac53}(\mathbb R^5)} + \|z\|_{L^{\frac52}(\mathbb R^5)} \le C\varepsilon^2.
\end{equation}

To recover the pressure, let \(R_i\) be the Riesz transform with Fourier symbol \(i\xi_i |\xi|^{-1}\), and define
\begin{align*}
T_z = \chi w\otimes w +\widetilde U\otimes z +z\otimes\widetilde U,
\qquad
q_z=\sum_{i,j=1}^{5}R_iR_j(T_z)_{ij}.
\end{align*}
For almost every \(s\in\mathbb R\), H\"older's inequality and Sobolev inequalities give
\begin{align*}
\|\widetilde U(s,\cdot)\otimes z(s,\cdot)\|_{L^{\frac53}(\mathbb R^4)} &\le\|\widetilde U(s,\cdot)\|_{L^4(\mathbb R^4)} \|z(s,\cdot)\|_{L^{\frac{20}{7}}(\mathbb R^4)}\\
&\le C\|\widetilde U(s,\cdot)\|_{L^4(\mathbb R^4)} \|\nabla_yz(s,\cdot)\|_{L^{\frac53}(\mathbb R^4)}.
\end{align*}
Taking the
\(L^{\frac53}(\mathbb R)\) norm in \(s\), we obtain
\begin{align*}
\|\widetilde U\otimes z\|_{L^{\frac53}(\mathbb R^5)}
+
\|z\otimes\widetilde U\|_{L^{\frac53}(\mathbb R^5)}
&\le
C\|\widetilde U\|_
 {L^\infty(\mathbb R;L^4(\mathbb R^4))}
\|\nabla_y z\|_{L^{\frac53}(\mathbb R^5)}
 \le C\delta\varepsilon^2.
\end{align*}
The \(L^{\frac53}(\mathbb R^5)\)-boundedness of the Riesz
transforms therefore implies
\begin{equation}\label{zpressure}
\|q_z\|_{L^{\frac53}(\mathbb R^5)}
\le C\|T_z\|_{L^{\frac53}(\mathbb R^5)}
\le C\varepsilon^2.
\end{equation}

The defining equation for \(z\) is
\begin{align*}
z=- S_a  T_z, \qquad\text{hence}\qquad L_a  z=-\mathbb P \nabla \cdot  T_z.
\end{align*}
Our choice of \(q_z\) satisfies
\begin{align*}
\Delta q_z = -\sum_{i,j=1}^{5}\partial_i\partial_j(T_z)_{ij}, \qquad \nabla q_z = -(I-\mathbb P )\nabla \cdot  T_z.
\end{align*}
Consequently,
\begin{equation}\label{zEquation}
L_a  z +\nabla \cdot (\widetilde U\otimes z+z\otimes\widetilde U) +\nabla q_z = -\nabla \cdot (\chi w\otimes w), 
\qquad 
\nabla \cdot  z=0 \quad\text{in }\mathbb R^5.
\end{equation}

{\underline{\bf Step 2: qualitative regularity of the linearized part.}} Let \(p_U\) denote the pressure associated with the smooth solution \(ae_1+U\). Subtracting its equation from the equation for \(u\), we find
\begin{align*}
 L_a w +\nabla \cdot (U\otimes w+w\otimes U) +\nabla(p-p_U) = -\nabla \cdot (w\otimes w),
\qquad
\nabla \cdot  w=0 \quad\text{in }Q_1.
\end{align*}
Set
\begin{align*}
h=w-z,
\qquad
\pi=p-p_U-q_z.
\end{align*}
Since \(\chi=1\) and \(\widetilde U=U\) on \({Q_\frac78}\),
subtracting \eqref{zEquation} gives
\begin{equation}\label{hEquation}
 L_a  h+U\cdot\nabla h+h\cdot\nabla U+\nabla\pi=0,
\qquad
\nabla \cdot  h=0
\quad\text{in }{Q_\frac78}.
\end{equation}
Interior Stokes regularity and bootstrapping applied to \eqref{hEquation}
give \(h,\pi\in C^\infty(Q_{\frac78})\).
Consequently,
\[
z=w-h\in H^1_{\mathrm{loc}}(Q_{\frac78}).
\]

{\underline{\bf Step 3: relative energy with the smooth comparison field.}} 
Define
\begin{align*}
V=ae_1+U+h, \qquad P=p_U+\pi.
\end{align*}
Adding the equation for \(ae_1+U\) to \eqref{hEquation}, we obtain
\begin{align*}
-\Delta V+(V\cdot\nabla)V+\nabla P = (h\cdot\nabla)h,
\qquad
\nabla \cdot  V=0 \quad\text{in }{Q_\frac78}.
\end{align*}
Moreover,
\begin{align*}
u-V=z, \qquad p-P=q_z.
\end{align*}
Since
\begin{align*}
(z\cdot\nabla V)\cdot z + ((h\cdot\nabla)h)\cdot z = (z\cdot\nabla U)\cdot z + (w\cdot\nabla h)\cdot z,
\end{align*}
we use Lemma \ref{relativeLemma} with \(F=(h\cdot\nabla)h\) to obtain, for every nonnegative \(\phi\in C_c^\infty({Q_\frac78})\),
\begin{align}
\int_{{Q_\frac78}}|\nabla z|^2\phi\,\mathrm dx \le{}& \frac12\int_{{Q_\frac78}}|z|^2\Delta\phi\,\mathrm dx + \frac a2\int_{{Q_\frac78}}|z|^2\partial_s\phi\,\mathrm dx
\notag\\
&+ \frac12\int_{{Q_\frac78}}|z|^2(U+w)\cdot\nabla\phi\,\mathrm dx + \int_{{Q_\frac78}} q_z z\cdot\nabla\phi\,\mathrm dx
\notag\\
&- \int_{{Q_\frac78}}(z\cdot\nabla U)\cdot z\,\phi\,\mathrm dx - \int_{{Q_\frac78}}(w\cdot\nabla h)\cdot z\,\phi\,\mathrm dx.
\label{relative}
\end{align}
%

{\underline{\bf Step 4: Preliminary energy estimate.}} Let
\begin{align*}
\phi=\psi^2, \qquad \psi\in C_c^\infty({Q_\frac78}), \qquad 0\le\psi\le1,
\end{align*}
in \eqref{relative} and rewrite the inequality as follows:
\begin{align*}
\int_{{Q_\frac78}}|\nabla z|^2 \psi^2 \,\mathrm dx \le{}& \frac12\int_{{Q_\frac78}}|z|^2\Delta \psi^2 \,\mathrm dx + \frac a2\int_{{Q_\frac78}}|z|^2\partial_s\psi^2 \,\mathrm dx
\notag\\
&+ \frac12\int_{{Q_\frac78}}|z|^2 U \cdot\nabla\psi^2 \,\mathrm dx + \frac12\int_{{Q_\frac78}}|z|^2 w \cdot\nabla\psi^2 \,\mathrm dx + \int_{{Q_\frac78}} q_z z\cdot\nabla \psi^2 \,\mathrm dx
\notag\\
&- \int_{{Q_\frac78}}(z\cdot\nabla U)\cdot z\,\psi^2 \,\mathrm dx - \int_{{Q_\frac78}}(w\cdot\nabla h)\cdot z\,\psi^2 \,\mathrm dx\\
&:= \sum_{k=1}^7I_k.
\end{align*}
We estimate the terms on the right-hand side.

{\underline{\it Estimate of $I_1$:}} By \eqref{weakz}, there holds
\begin{align*}
I_1 \le \|\Delta\phi\|_{L^\infty({Q_\frac78})} \|z\|_{L^2(\mathbb R^5)}^2 \le Ca^{-1}\|\Delta\phi\|_{L^\infty({Q_\frac78})} \varepsilon^4.
\end{align*}

{\underline{\it Estimate of $I_2$:}} Since $\partial_s \psi^2 = 2\psi \partial_s \psi$, using \eqref{weakz} again, there holds
\begin{align*}
I_2 \le a\|\partial_s\phi\|_{L^\infty(Q_\frac78)} \|z\|_{L^2(\mathbb R^5)}^2 \le C\|\partial_s\phi\|_{L^\infty(Q_\frac78)} \varepsilon^4.
\end{align*}

{\underline{\it Estimate of $I_3$:}} Since $\psi$ has compact support, $z\psi$ can be extended to the whole space. Using the Sobolev inequality, one obtains
\begin{align*}
I_3 \leq \|U\|_{L^\infty((-\frac78,\frac78);L^4(B_{\frac78}^4))} \|z\psi\|_{L^2(\mathbb R;L^4(\mathbb R^4))} \|z\nabla\psi\|_{L^2(Q_\frac78)}.
\end{align*}
The homogeneous Sobolev inequality gives
\begin{align}\label{ine:zpsi24} \nonumber
\|z\psi\|_{L^2(\mathbb R;L^4(\mathbb R^4))}^2 &= \int_{\mathbb R} \|z(s,\cdot)\psi(s,\cdot)\|_{L^4(\mathbb R^4)}^2 \,\mathrm ds \le C\int_{\mathbb R}
\|\nabla_y(z\psi)(s,\cdot)\|_{L^2(\mathbb R^4)}^2 \,\mathrm ds\\
&= C\|\nabla_y(z\psi)\|_{L^2(\mathbb R^5)}^2 \le C\left(\|\nabla z \psi\|_{L^2(Q_\frac78)}^2 + \|z\nabla\psi\|_{L^2(Q_\frac78)}^2\right).
\end{align}
Also, 
\ben\label{ine:znablapsi2}
\|z\nabla\psi\|_{L^2(Q_\frac78)}^2 \le \|\nabla\psi\|_{L^\infty(Q_\frac78)}^2 \|z\|_{L^2(\mathbb R^5)}^2 \le Ca^{-1}\|\nabla\psi\|_{L^\infty(Q_\frac78)}^2
\varepsilon^4.
\een
\eqref{ass:keyprop}, \eqref{ine:zpsi24}, \eqref{ine:znablapsi2} and Young's inequality yield
\begin{align*}
I_3 &\leq C \delta \left(\|\nabla z \psi\|_{L^2(Q_\frac78)}^2 + a^{-1}\|\nabla\psi\|_{L^\infty(Q_{\frac78})}^2\varepsilon^4\right)^\frac12 a^{-\frac12} \|\nabla \psi\|_{L^\infty(Q_\frac78)} \varepsilon^2 \\
&\leq \frac 1{16}\|\nabla z \psi\|_{L^2(Q_\frac78)}^2 + C (\delta + \delta^2) a^{-1}\|\nabla\psi\|_{L^\infty(Q_\frac78)}^2 \varepsilon^4.
\end{align*}

{\underline{\it Estimate of $I_5$:}} By \eqref{weakz} and \eqref{zpressure}, we have
\begin{align*}
I_5 \leq \|\nabla\phi\|_{L^\infty(Q_\frac78)} \|q_z\|_{L^{\frac53}(\mathbb R^5)} \|z\|_{L^{\frac52}(\mathbb R^5)} \le C\|\nabla\phi\|_{L^\infty(Q_\frac78)}
\varepsilon^4.
\end{align*}

{\underline{\it Estimate of $I_6$:}} By \eqref{ass:keyprop}, \eqref{ine:znablapsi2} and \eqref{ine:zpsi24}, we have
\begin{align*}
I_6 &\leq \int_{{Q_\frac78}}|\nabla U|\,|z\psi|^2\,\mathrm dx \le \|\nabla U\|_ {L^\infty((-\frac78,\frac78);L^2(B_{\frac78}^4))} \|z\psi\|_ {L^2(\mathbb R;L^4(\mathbb R^4))}^2 \\
&\leq C \delta \left(\|\nabla z \psi\|_{L^2(Q_\frac78)}^2 + \|z\nabla\psi\|_{L^2(Q_\frac78)}^2\right) \\
&\leq C \delta \|\nabla z \psi\|_{L^2(Q_\frac78)}^2 + C \delta a^{-1}\|\nabla\psi\|_{L^\infty(Q_\frac78)}^2 \varepsilon^4.
\end{align*}

{\underline{\it Estimate of $I_4+I_7$:}} Since the terms $I_4$ and $I_7$ contain $w$, we treat these two terms together. Letting $z=w-h$, we expand the two terms containing \(w\), where $h$ is smooth on the support of \(\phi\). Since $\nabla \cdot w = 0$, there holds
\begin{align*}
\int_{{Q_\frac78}}(w\cdot\nabla h)\cdot h\,\phi\,\mathrm dx = \frac12\int_{{Q_\frac78}}w\cdot\nabla(|h|^2)\phi\,\mathrm dx = -\frac12\int_{{Q_\frac78}}|h|^2w\cdot\nabla\phi\,\mathrm dx.
\end{align*}
Using \(w\cdot h=|w|^2-w\cdot z\), we obtain 
\begin{align}\nonumber
I_4 + I_7&= -\frac12\int_{{Q_\frac78}}|w|^2w\cdot\nabla\phi\,\mathrm dx + \int_{{Q_\frac78}}(w\cdot z)w\cdot\nabla\phi\,\mathrm dx - \int_{{Q_\frac78}}(w\cdot\nabla h)\cdot w\,\phi\,\mathrm dx \\
&:= J_1 + J_2 + J_3.
\label{firstidentity}
\end{align}
By \eqref{ine:w103}, \eqref{weakz} and the fact that $\phi$ has compact support, one obtains
\begin{align*}
J_1 + J_2 &\leq \|\nabla\phi\|_{L^\infty(Q_\frac78)} \|w\|_{L^3(Q_1)}^3 + \|\nabla\phi\|_{L^\infty(Q_\frac78)} \|w\|_{L^{\frac{10}{3}}(Q_1)}^2 \|z\|_{L^{\frac52}(\mathbb R^5)} \\
&\leq C\|\nabla\phi\|_{L^\infty(Q_\frac78)}\varepsilon^3 + C\|\nabla\phi\|_{L^\infty(Q_\frac78)}\varepsilon^4.
\end{align*}
Next, we deal with the term $J_3$. Let \(\frac23 \le r<R\le\frac34\). Choose \(\psi\) such that
\begin{align*}
0\le\psi\le1, \qquad \psi=1\ \text{on }Q_r, \qquad \psi\in C_c^\infty(Q_{(r+R)/2}),
\end{align*}
and
\begin{align*}
\|\nabla\psi\|_{L^\infty(Q_\frac78)} \le C(R-r)^{-1}, \qquad \|\nabla^2\psi\|_{L^\infty(Q_\frac78)} \le C(R-r)^{-2}.
\end{align*}
Taking \(\delta\) below the threshold in Lemma \ref{linearized}, we may apply that lemma to \eqref{hEquation} between \(Q_R\) and \(Q_{\frac{r+R}2}\). It yields
\begin{align}\label{ine:yahq}
\Y_a(h;Q_{\frac{r+R}2})\le C(R-r)^{-M} \|h\|_{H^1(Q_R)}.
\end{align}
By \eqref{ine:wL2L4} and \eqref{ine:yahq}, we have
\begin{align*}
J_3&\le \int_{-\frac{r+R}2}^{\frac{r+R}2} \|\nabla h(s,\cdot)\|_{L^2(B_{\frac{r+R}2}^4)} \|w(s,\cdot)\|_{L^4(B_{\frac{r+R}2}^4)}^2 \,\mathrm ds\\
&\le \|\nabla h\|_
 {L^\infty((-\frac{r+R}2,\frac{r+R}2); L^2(B_{\frac{r+R}2}^4))} \|w\|_{L^2((-1,1);L^4(B_1^4))}^2\\
&\le C(R-r)^{-M}\varepsilon^2 \|h\|_{H^1(Q_R)}.
\end{align*}
Since $h = -z + w$, by \eqref{ass:keyprop} and \eqref{weakz}, we have 
\begin{align*}
\|h\|_{H^1(Q_R)} \le \|w\|_{H^1(Q_R)} + \|z\|_{L^2(Q_R)} + \|\nabla z\|_{L^2(Q_R)} \le \varepsilon + Ca^{-\frac12}\varepsilon^2+ \|\nabla z\|_{L^2(Q_R)},
\end{align*}
which means that 
\begin{align*}
J_3 \leq C(R-r)^{-M}\varepsilon^2 \left(\varepsilon + Ca^{-\frac12}\varepsilon^2+ \|\nabla z\|_{L^2(Q_R)}\right).
\end{align*}
Collecting the estimates of $J_1$, $J_2$ and $J_3$, we have 
\begin{align*}
I_4 + I_7 &\leq C\|\nabla\phi\|_{L^\infty(Q_\frac78)}\varepsilon^3 + C\|\nabla\phi\|_{L^\infty(Q_\frac78)}\varepsilon^4 \\
&+ C(R-r)^{-M}\varepsilon^2 \left(\varepsilon + Ca^{-\frac12}\varepsilon^2+ \|\nabla z\|_{L^2(Q_R)}\right).
\end{align*}

Collecting the estimates of $I_1$ to $I_7$ and noting that $C \delta < \frac14$, we obtain
\begin{align*}
\|\nabla z\|_{L^2(Q_r)}^2 \leq \frac12 \|\nabla z\|_{L^2(Q_R)}^2  +  C (R-r)^{-M} \varepsilon^2 \left(\varepsilon + \|\nabla z\|_{L^2(Q_R)}\right).
\end{align*}
Using Young's inequality, there holds
\begin{align*}
\|\nabla z\|_{L^2(Q_r)}^2 \leq \frac34 \|\nabla z\|_{L^2(Q_R)}^2 +   C (R-r)^{-M} \varepsilon^3 +  C (R-r)^{-2M} \varepsilon^4.
\end{align*}
Since $z\in H^1_{\mathrm{loc}}(Q_{\frac78})$, applying Lemma~\ref{lem:Gia} on \([\tfrac23,\tfrac34]\), we obtain
\begin{align*}
\|\nabla z\|_{L^2(Q_\frac23)}^2 \leq  C\varepsilon^3.
\end{align*}
Combining this estimate with \eqref{weakz}, we find
\ben\label{zfirstpass}
\|z\|_{H^1(Q_\frac23)}^2 \leq Ca^{-1}\varepsilon^4+  C\varepsilon^3.
\een
Using $h=w-z$, \eqref{ass:keyprop}, \eqref{zfirstpass} and Lemma \ref{linearized}, we obtain
\ben\label{firstpass}
\|h\|_{H^1(Q_\frac23)} \leq C\varepsilon^\frac32 + C \varepsilon \leq C \varepsilon, \qquad \Y_a(h;Q_\frac12) \leq C \varepsilon.
\een

{\underline{\bf Step 5: the second energy estimate.}}  After obtaining the $H^1$ norm of $h$ from \eqref{firstpass}, the estimate for $I_4 + I_7$ can be improved.

Choose a fixed cutoff
\begin{align*}
\psi\in C_c^\infty(Q_{\frac12}),
\qquad
0\le\psi\le1,
\qquad
\psi=1\ \text{on }Q_{\frac38},
\end{align*}
and again set \(\phi=\psi^2\). All derivatives of this cutoff are bounded by universal constants. Moreover, the estimates for $I_1,I_2,I_3,I_5,I_6$ follow from the preceding step.

In \eqref{relative}, write \(w=h+z\). The terms that were treated by \eqref{firstidentity} now become
\begin{align*}
&\frac12\int_{{Q_\frac78}}|z|^2w\cdot\nabla\phi\,\mathrm dx - \int_{{Q_\frac78}}(w\cdot\nabla h)\cdot z\,\phi\,\mathrm dx \\
&\qquad= \int_{{Q_\frac78}}|z|^2h\cdot\psi\nabla\psi\,\mathrm dx + \frac12\int_{{Q_\frac78}}|z|^2z\cdot\nabla\phi\,\mathrm dx \\
&\qquad\quad - \int_{{Q_\frac78}}(z\cdot\nabla h)\cdot z\,\phi\,\mathrm dx - \int_{{Q_\frac78}}(h\cdot\nabla h)\cdot z\,\phi\,\mathrm dx \\
&\quad := K_1 + K_2 + K_3 + K_4.
\end{align*}

{\underline{\it Estimate of $K_1$:}} By \eqref{ine:zpsi24}, \eqref{ine:znablapsi2} and \eqref{firstpass}, there holds
\begin{align*}
K_1 &\leq \|h\|_ {L^\infty((-\frac12,\frac12);L^4(B_{\frac12}^4))}\|z\psi\|_ {L^2(\mathbb R;L^4(\mathbb R^4))}\|z\nabla \psi\|_ {L^2(\mathbb R;L^2(\mathbb R^4))} \\
&\leq C \varepsilon^3 \left(\|\nabla z \psi\|_{L^2(Q_\frac12)}^2 + \varepsilon^4\right)^\frac12 \\
&\leq \frac1{16} \|\nabla z \psi\|_{L^2(Q_\frac12)}^2 + C \varepsilon^6 + C \varepsilon^5.
\end{align*}

{\underline{\it Estimate of $K_2$:}} By H\"older's inequality, Sobolev inequality, \eqref{weakz} and \eqref{zfirstpass}, we obtain
\begin{align*}
K_2 &\le C\int_{Q_{\frac12}}|z|^3\,\mathrm dx \le C\|z\|_{L^{\frac52}(Q_{\frac12})} \|z\|_{L^{\frac{10}{3}}(Q_{\frac12})}^2\\
&\le C\|z\|_{L^{\frac52}(\mathbb R^5)}\|z\|_{H^1(Q_{\frac23})}^2 \le C\varepsilon^5.
\end{align*}

{\underline{\it Estimate of $K_3$:}} By Young's inequality, \eqref{ine:zpsi24}, \eqref{ine:znablapsi2}, \eqref{firstpass} and noting that $C \varepsilon < \frac1{16}$, there holds
\begin{align*}
K_3 &\le \|\nabla h\|_ {L^\infty((-\frac12,\frac12);L^2(B_{\frac12}^4))}\|z\psi\|_ {L^2(\mathbb R;L^4(\mathbb R^4))}^2 \le C\varepsilon \left(\|\nabla z \psi\|_{L^2(Q_\frac12)}^2 + \varepsilon^4\right) \\
&\leq \frac1{16} \|\nabla z \psi\|_{L^2(Q_\frac12)}^2 + C\varepsilon^5.
\end{align*}

{\underline{\it Estimate of $K_4$:}} The ordinary \(L^4(Q_{\frac12})\) norm of \(h\) satisfies
\begin{align}\label{ine:step5esth}
\|h\|_{L^4(Q_{\frac12})}^4 &= \int_{-\frac12}^{\frac12} \|h(s,\cdot)\|_{L^4(B_{\frac12}^4)}^4 \,\mathrm ds  \le \|h\|_ {L^\infty((-\frac12,\frac12);L^4(B_{\frac12}^4))}^4\le C\varepsilon^4.
\end{align}
By integration by parts, one obtains
\begin{align*}
-\int_{{Q_\frac78}}(h\cdot\nabla h)\cdot z\,\phi\,\mathrm dx &=-\int_{{Q_\frac78}} h_j\partial_j h_i\,z_i\phi\,\mathrm dx\\
&=\int_{{Q_\frac78}} h_i h_j\partial_j z_i\,\phi\,\mathrm dx+\int_{{Q_\frac78}} h_i h_j z_i\partial_j\phi\,\mathrm dx.
\end{align*}
Then 
\begin{align*}
K_4 \leq C \|h\|_{L^4(Q_{\frac12})}^2 \|\nabla z\psi\|_{L^2(Q_\frac12)} +  C \|h\|_{L^4(Q_{\frac12})}^2 \|z\|_{L^2(Q_\frac12)}.
\end{align*}
Young's inequality, \eqref{weakz} and \eqref{ine:step5esth} yield
\begin{align*}
K_4 \leq \frac1{16} \|\nabla z \psi\|_{L^2(Q_\frac12)}^2 + C \varepsilon^4.
\end{align*}
Combining the estimates for $K_1$ to $K_4$ and recalling the estimate for $I_1,I_2,I_3,I_5,I_6$, we obtain 
\begin{equation}\label{secondpass}
\|z\|_{H^1(Q_{\frac38})} \le C\varepsilon^2.
\end{equation}

{\underline{\bf Step 6: correction to an exact smooth background.}} On \(Q_{\frac12}\), the field \(ae_1+U\) is an exact smooth solution, while \(h\) is a smooth divergence-free solution of \eqref{hEquation}.
Moreover,
\begin{align*}
\Y_a(U;Q_{\frac12})\le\delta,
\qquad
\Y_a(h;Q_{\frac12})\le C_0\varepsilon
\end{align*}
for a universal constant \(C_0\), by \eqref{firstpass}. For sufficiently small \(\bar\delta\) and \(\bar\varepsilon\),
Lemma~\ref{strongcorrect} gives a divergence-free corrector \(k\) such that \(\Y_a(k)\le C\varepsilon^2\) and \(ae_1+U+h+k\) is an exact smooth solution in \(Q_{\frac14}\).

Set \(U_+=U+h+k\). Then
\begin{align*}
u-ae_1-U_+=w-h-k=z-k \quad\text{on }Q_{\frac14}.
\end{align*}
The local \(H^1\) norm of \(k\) is controlled directly by its first two mixed norms:
\begin{align*}
\|k\|_{H^1(Q_{\frac14})}^2 &= \int_{-\frac14}^{\frac14} \int_{B_{\frac14}^4} \bigl(|k(s,y)|^2+|\nabla k(s,y)|^2\bigr)\,\mathrm dy\,\mathrm ds\\
&\le \frac12|B_{\frac14}^4|^{\frac12} \|k\|_{L^\infty(\mathbb R;L^4(\mathbb R^4))}^2 + \frac12 \|\nabla k\|_ {L^\infty(\mathbb R;L^2(\mathbb R^4))}^2 \le
C\Y_a(k)^2 \le C\varepsilon^4.
\end{align*}
In particular, this estimate introduces no factor depending on \(a\). Together with \eqref{secondpass}, it yields
\begin{align*}
\|u-ae_1-U_+\|_{H^1(Q_{\frac14})} \le \|z\|_{H^1(Q_{\frac14})} + \|k\|_{H^1(Q_{\frac14})} \le C\varepsilon^2.
\end{align*}
Finally,
\begin{align*}
\Y_a(U_+;Q_{\frac14}) &\le \Y_a(U;Q_{\frac14}) + \Y_a(h;Q_{\frac14}) + \Y_a(k;Q_{\frac14})  \le \delta+C\varepsilon+C\varepsilon^2 \le \delta+C'\varepsilon.
\end{align*}
\end{proof}

\section{Proofs of the main results}
\subsection{The one-scale regularity criterion}
\begin{proof}[Proof of Theorem~\ref{thm:main}]
Let \(\bar\delta,\bar\varepsilon>0\) be the thresholds in Proposition~\ref{improvement}, and let \(C_{\mathrm{imp}}\ge1\) be a universal constant that bounds both constants in its conclusion. The energy threshold \(\varepsilon_\ast\) will be fixed below.

{\underline{\bf Step 1: normalization and the initial decomposition.}} 
Define
\begin{align*}
v(y)=R u(Ry), \qquad q(y)=R^2p(Ry), \qquad y\in B_1.
\end{align*}
By scaling, $(v,q)$ is a suitable weak solution in $B_1$, and
\begin{align*}
\int_{B_1}|\nabla v(y)|^2\,\mathrm dy=\frac1R\int_{B_R}|\nabla u(x)|^2\,\mathrm dx := E_u(R).
\end{align*}
For the constant $b=(v)_{B_1}$, by Poincar\'e's inequality, we obtain
\begin{equation}\label{initialPoincare}
\|v-b\|_{H^1(B_1)} \le C\|\nabla v\|_{L^2(B_1)} = C E_u(R)^{\frac12}.
\end{equation}
By rotational invariance, we may assume \(b=|b|e_1\).
Let $u_0(x) = \frac12 v(\frac x2)$ and $p_0(x)=\frac14 q(\frac x2)$ for $x\in Q_1$ (so $\frac x2\in B_1$), and set
\begin{align*}
a_0 = \frac{|b|}2, \qquad U_0 = 0, \qquad w_0 = u_0 - a_0 e_1,
\end{align*}
so that $u_0 = a_0 e_1 + U_0 + w_0$ and $a_0 e_1 + U_0$ is a smooth solution of $\eqref{NS}$. \eqref{initialPoincare} implies
\begin{align*}
\|w_0\|_{H^1(Q_1)}^2 &\leq C \int_{Q_{\frac12}}|v(y)-b|^2\,\mathrm dy + C \int_{Q_{\frac12}}|\nabla v(y)|^2\,\mathrm dy\\
&\le C\|v-b\|_{H^1(B_1)}^2 \le C E_u(R).
\end{align*}
Now for $\delta_0 = 0$ and $\varepsilon_0 = C_{\mathrm{init}}\varepsilon_\ast$ with a fixed universal constant $C_{\mathrm{init}}>1$, we have
\begin{equation}\label{initialBounds}
\|U_0\|_{L^\infty((-1,1);L^4(B_1^4))} \leq \delta_0, \qquad \|w_0\|_{H^1(Q_1)} \le \varepsilon_0.
\end{equation}

{\underline{\bf Step 2: mathematical induction.}} Fix \(\theta=\frac14\) and let \(N=\min\{j\in\mathbb{N}_0:\theta^j a_0\le1\}\). If \(N=0\), proceed to Step~3.
For \(0\le j<N\), assume inductively that we have a suitable weak solution
\begin{align*}
u_j=a_je_1+U_j+w_j \quad\text{in }Q_1,
\end{align*}
where \(a_je_1+U_j\) is a smooth solution of $\eqref{NS}$ and
\begin{align*}
\Y_{a_j}(U_j;Q_1)\le\delta_j\le\bar\delta, \qquad \|w_j\|_{H^1(Q_1)} \le\varepsilon_j\le\bar\varepsilon.
\end{align*}
Since $a_j=\theta^j a_0>1$, Proposition \ref{improvement} gives a smooth solution \(a_je_1+U_{j,+}\) on \(Q_\theta\).
Writing \(w_{j,+}=u_j-a_je_1-U_{j,+},\) we have
\begin{align}\label{ine:mathinduj}
\Y_{a_j}(U_{j,+};Q_\theta) \le \delta_j+C_{\mathrm{imp}}\varepsilon_j,
\qquad
\|w_{j,+}\|_{H^1(Q_\theta)} \le C_{\mathrm{imp}}\varepsilon_j^2.
\end{align}
With $C_{\mathrm{imp}}$ fixed as above, rescale \(Q_\theta\) to \(Q_1\) by setting
\begin{align*}
u_{j+1}(x)&=\theta u_j(\theta x),
&
p_{j+1}(x)&=\theta^2p_j(\theta x),
\\
U_{j+1}(x)&=\theta U_{j,+}(\theta x),
&
w_{j+1}(x)&=\theta w_{j,+}(\theta x).
\end{align*}
Then
\begin{align*}
u_{j+1}=a_{j+1}e_1+U_{j+1}+w_{j+1},\qquad a_{j+1}=\theta a_j.
\end{align*}
By \eqref{Yscale} and \eqref{ine:mathinduj}, there holds
\begin{align*}
\|w_{j+1}\|_{H^1(Q_1)} &\le\theta^{-\frac32}\|w_{j,+}\|_{H^1(Q_\theta)} \le C_{\mathrm{imp}}\theta^{-\frac32} \varepsilon_j^2,\\ 
\Y_{a_{j+1}}(U_{j+1};Q_1) &=\Y_{a_j}(U_{j,+};Q_\theta)
 \le\delta_j+C_{\mathrm{imp}}\varepsilon_j.
\end{align*}
Set
\begin{equation}\label{iteration}
a_{j+1}=\theta a_j, \qquad \varepsilon_{j+1} =C_{\mathrm{imp}}\theta^{-\frac32}\varepsilon_j^2, \qquad \delta_{j+1} = \delta_j+C_{\mathrm{imp}}\varepsilon_j,
\end{equation}
which, provided the bounds verified below hold, gives by induction a suitable weak solution at each step $1\le n\le N$ satisfying
\begin{align*}
u_n=a_ne_1+U_n+w_n \quad\text{in }Q_1,
\end{align*}
where \(a_ne_1+U_n\) is a smooth solution of $\eqref{NS}$ and
\begin{align*}
\Y_{a_n}(U_n;Q_1)\le\delta_n\le\bar\delta, \qquad \|w_n\|_{H^1(Q_1)} \le\varepsilon_n\le\bar\varepsilon.
\end{align*}
To verify these bounds, assume that
\begin{align}\label{ine:cimpbdd}
\varepsilon_0\le\bar\varepsilon, \qquad C_{\mathrm{imp}} \leq \min\{\frac12 \theta^\frac32 \varepsilon_0^{-1}, \frac14 \bar\delta \varepsilon_0^{-1}\}\NSrev{.}
\end{align}
Then \eqref{iteration} and \eqref{ine:cimpbdd} give
\begin{align*}
\varepsilon_j = (C_{\mathrm{imp}}\theta^{-\frac32})^{2^j-1} \varepsilon_0^{2^j}=(C_{\mathrm{imp}}\theta^{-\frac32})^{-1}(C_{\mathrm{imp}}\theta^{-\frac32}\varepsilon_0)^{2^j} \leq 2^{-2^j+1} \varepsilon_0.
\end{align*}
Since $-2^j+1 < -j$ for $j > 1$, we have 
\begin{align}\label{eq:sumepsilon0}
\sum_{j=0}^{N}\varepsilon_j \le2\varepsilon_0.
\end{align}
By \eqref{eq:sumepsilon0}, \eqref{iteration} and \eqref{ine:cimpbdd}, we have 
\begin{align*}
\delta_j=\delta_0+C_{\mathrm{imp}}\sum_{i=0}^{j-1}\varepsilon_i \le2C_{\mathrm{imp}}\varepsilon_0\le\frac12\bar\delta.
\end{align*}
Thus the bounds required at every $j<N$ hold independently of $N$. The scaling identity \eqref{Yscale} is valid also at the last step, when $a_N\le1$. We have
\begin{align}\label{ine:smalldeltavarepsilonN}
0\le a_N\le1,\qquad \delta_N\le2C_{\mathrm{imp}}\varepsilon_0,\qquad \varepsilon_N\le\varepsilon_0.
\end{align}

{\underline{\bf Step 3: smallness of the velocity.}} For \(u_N\), we have
\begin{align*}
u_N=a_Ne_1+U_N+w_N, \qquad 0\le a_N\le1,
\end{align*}
with
\begin{align*}
\|U_N\|_{L^\infty((-1,1);L^4(B_1^4))} \le\delta_N,
\qquad
\|w_N\|_{H^1(Q_1)} \le\varepsilon_N.
\end{align*}
Let \(0<\rho<\frac14\). Next, we compute the $L^3$ norm of $u_N$ on the ball $B_\rho$ at one scale.

For the part $a_N e_1$, since $a_N$ is a constant, there holds
\begin{align}\label{ine:finalan}
\rho^{-2}\|a_Ne_1\|_{L^3(B_\rho)}^3 = \rho^{-2}|a_N|^3|B_\rho| = |B_1|\rho^3 |a_N|^3 \le C\rho^3.
\end{align}

For the part $U_N$, $B_\rho\subset Q_\rho$ and H\"older's inequality gives
\begin{align}\label{ine:finalun}\nonumber
\rho^{-2}\int_{B_\rho}|U_N(x)|^3\,\mathrm dx &\le \rho^{-2} \int_{-\rho}^{\rho} \int_{B_\rho^4}|U_N(s,y)|^3 \,\mathrm dy\,\mathrm ds \\ \nonumber
&\le \rho^{-2} |B_\rho^4|^{\frac14} \int_{-\rho}^{\rho} \|U_N(s,\cdot)\|_{L^4(B_\rho^4)}^3 \,\mathrm ds\\
&\le 2\rho^{-1}\,|B_\rho^4|^{\frac14} \|U_N\|_{L^\infty((-1,1);L^4(B_1^4))}^3 \le C\delta_N^3.
\end{align}

For the part $w_N$, the Sobolev embedding gives
\begin{align}\label{ine:finalwn}
\rho^{-2}\|w_N\|_{L^3(B_\rho)}^3 &\le \rho^{-2}|B_\rho|^{\frac{1}{10}} \|w_N\|_{L^{\frac{10}{3}}(B_\rho)}^3 \le C\rho^{-\frac32} \|w_N\|_{L^{\frac{10}{3}}(Q_1)}^3 \le C\rho^{-\frac32}\varepsilon_N^3.
\end{align}

Combining \eqref{ine:finalan}, \eqref{ine:finalun} and \eqref{ine:finalwn}, fix a universal constant \(C_{\mathrm{vel}}\ge1\) such that
\begin{equation}\label{terminal}
\rho^{-2}\|u_N\|_{L^3(B_\rho)}^3 \le C_{\mathrm{vel}} \bigl(\rho+\delta_N+\rho^{-\frac12}\varepsilon_N \bigr)^3.
\end{equation}

For the constant $\varepsilon_{\mathrm{cui}}$ in Theorem~1.1 of \cite{Cui2023} (with $q=3$ and $f=0$), fix $\rho_\ast \in (0,\frac14)$ satisfying
\begin{align*}
C_{\mathrm{vel}} \rho_*^3 \leq \frac14 \varepsilon_{\mathrm{cui}}.
\end{align*}
Choose $\varepsilon_*>0$ so small that $\varepsilon_0=C_{\rm init}\varepsilon_*$ satisfies
\begin{equation}\label{finalThreshold}
\varepsilon_0< \min\left\{ \bar\varepsilon,\, \frac{\theta^{\frac32}}{2C_{\mathrm{imp}}},\, \frac{\bar\delta}{4C_{\mathrm{imp}}},\, \frac{\varepsilon_{\mathrm{cui}}^{\frac13}} {4C_{\mathrm{vel}} (2C_{\mathrm{imp}}+\rho_*^{-\frac12})}
\right\}\NSrev{.}
\end{equation}
Then \eqref{terminal} gives
\begin{align*}
\rho_*^{-2} \int_{B_{\rho_*}} |u_N|^3 dx <\varepsilon_{\mathrm{cui}}.
\end{align*}
After rescaling \(B_{\rho_*}\) to \(B_1\), Theorem~1.1 of
\cite{Cui2023} with \(q=3\) yields regularity of \(u_N\) at the origin, and hence of \(u\) at \(x_0\).

\end{proof}

\subsection{Minkowski dimension and content of the singular set}

\begin{proof}[The proof of Corollary~\ref{minkowski}]
For a nonempty bounded set \(S\subset\mathbb R^5\), let \(N(S,r)\) be the smallest number of balls of radius \(r\) covering \(S\). Its upper Minkowski dimension is
\begin{align*}
\overline{\dim}_{\mathrm M}S=\limsup_{r\downarrow0}\frac{\log N(S,r)}{\log(1/r)}.
\end{align*}
Setting \(S_r=\{x:\operatorname{dist}(x,S)<r\}\), with the normalization in Corollary~\ref{minkowski} we have
\begin{align*}
\mathcal M^1(S)=0 \quad\Longleftrightarrow\quad \lim_{r\downarrow0}r^{-4}|S_r|=0.
\end{align*}
Fix a compact set \(K\Subset\Omega\) and choose \(r_0>0\) so that
\begin{align*}
\overline{K_{4r_0}}\Subset\Omega, \qquad K_t=\{x\in\R^5:\operatorname{dist}(x,K)<t\}.
\end{align*}
If \(S_K=\varnothing\), the conclusions are immediate. Otherwise \(S_K\) is compact because the regular set is relatively open.
By Theorem \ref{thm:main}, if $x \in S_K$, then 
\begin{equation}\label{singularEnergyLowerBound}
\int_{B_r(x)} |\nabla u(y)|^2 \,dy \ge \varepsilon_\ast^2 r \qquad  0<r<r_0.
\end{equation}
Fix \(0<r<r_0\) and choose a maximal finite family \(x_1,\ldots,x_{M_r}\in S_K\) with \(|x_i-x_j|\ge2r\) for \(i\ne j\). The balls \(B_r(x_i)\) are
pairwise disjoint. Maximality and the triangle inequality imply
\begin{align*}
S_K\subset\bigcup_{i=1}^{M_r}B_{2r}(x_i),
\qquad
(S_K)_r\subset\bigcup_{i=1}^{M_r}B_{3r}(x_i),
\qquad
\bigcup_{i=1}^{M_r}B_r(x_i)\subset(S_K)_r.
\end{align*}
Summing \eqref{singularEnergyLowerBound} over the disjoint balls yields
\begin{equation}\label{packingEnergy}
\varepsilon_\ast^2 r M_r \le\sum_{i=1}^{M_r}\int_{B_r(x_i)} |\nabla u|^2 \, dx \le\int_{(S_K)_r} |\nabla u|^2  \, dx.
\end{equation}
Writing \(\omega_5=|B_1|\), we obtain
\begin{align}
&N(S_K,2r)\le M_r \le\frac{1}{\varepsilon_*^2r}\int_{(S_K)_r} |\nabla u|^2  \, dx, \label{coveringEstimate}\\
&|(S_K)_r|\le3^5\omega_5 r^5M_r \le\frac{3^5\omega_5}{\varepsilon_\ast^2}\,r^4 \int_{(S_K)_r} |\nabla u|^2  \, dx. \label{tubeProof}
\end{align}
\eqref{coveringEstimate} and \eqref{tubeProof} prove \eqref{tubeEstimate}. 
Since $(S_K)_r\subset K_{r_0}$ for $r<r_0$,
\eqref{coveringEstimate} gives $N(S_K,2r)\le \frac{C_K}{r}$,
and hence $\overline{\dim}_{\mathrm M}S_K\le1$.
Moreover, \eqref{tubeProof} implies
\(|(S_K)_r|\le 3^5\omega_5C_Kr^4\to0\).
The absolute continuity of the integral therefore yields
\begin{align*}
\int_{(S_K)_r} |\nabla u|^2 \,\mathrm dx\longrightarrow0.
\end{align*}
Dividing \eqref{tubeProof} by \(\omega_4r^4\), we obtain
\begin{align*}
0\le\mathcal M^1(S_K) \le\frac{3^5\omega_5}{\omega_4\varepsilon_\ast^2} \limsup_{r\downarrow0} \int_{(S_K)_r} |\nabla u|^2 \,\mathrm dx =0.
\end{align*}
This proves \eqref{MinkowskiConclusion}.
\end{proof}

\appendix 
\section{Auxiliary estimates}
\begin{lemma}\label{aux:H}
For the field \(H\) in the proof of Lemma~\ref{homogeneous}, we have the following estimate:
\begin{align*}
\|H\|_{H^2(Q_{\rho_6})}\le C(R-r)^{-M}\|h\|_{H^1(Q_R)}.
\end{align*}
\end{lemma}
\begin{proof}We now estimate the full $H^2$ norm of the auxiliary field $H$.
Recall that
\[
\Phi(s,y)
=
-a^{-1}\int_0^s\bigl(\pi(t,y)-\pi(0,0)\bigr)\,dt+\psi(y),
\qquad
H=\nabla\Phi.
\]
We write $H=(H_0,H_1,\ldots,H_4)$, where $H_0$ denotes the component in
the $s$-direction and the remaining components are in the $y$-directions.
Thus,
\[
H_0
=
-a^{-1}\bigl(\pi(s,y)-\pi(0,0)\bigr),
\]
and, for $\alpha\in\{1,\ldots,4\}$,
\[
H_\alpha
=
-a^{-1}\int_0^s\partial_{y_\alpha}\pi(t,y)\,dt
+
\partial_{y_\alpha}\psi(y).
\]

The pressure estimates above imply, after increasing the universal
exponent $M$ if necessary,
\[
\sum_{k=0}^{3}
{(R-r)}^k
\|D^k(\pi-\pi(0,0))\|_{L^\infty(Q_{\rho_7})}
\leq
Ca{(R-r)}^{-M}\|h\|_{H^1(Q_R)},
\]
where for $k=0$ the estimate is understood for the oscillation of $\pi$,
which is all that occurs in the definition of $\Phi$. In particular,
\[
\|D^k(\pi-\pi(0,0))\|_{L^\infty(Q_{\rho_7})}
\leq
Ca{(R-r)}^{-M}\|h\|_{H^1(Q_R)},
\qquad
0\leq k\leq 3.
\]

Recall that $\psi$ solves
\[
\Delta_y\psi
=
a^{-1}\partial_s\pi(0,\cdot)
\quad\text{in }B^4_{\rho_7}.
\]
Interior estimates for the Poisson equation, together with the pressure
bounds, give
\[
\sum_{j=1}^{3}
{(R-r)}^j
\|D_y^j\psi\|_{L^\infty(B^4_{\rho_6})}
\leq
C{(R-r)}^{-M}\|h\|_{H^1(Q_R)}.
\]
The factor $a^{-1}$ in the equation for $\psi$ cancels the factor $a$
in the pressure estimate. In particular, the constant in this estimate
is independent of $a\geq1$.

We first estimate $H$ itself. From the component formulas above, and
since $|s|\leq \rho_6\leq 1$, we obtain
\[
\begin{aligned}
|H_0(s,y)|
&\leq
a^{-1}
\left|\pi(s,y)-\pi(0,0)\right|
\leq
C{(R-r)}^{-M}\|h\|_{H^1(Q_R)},
\\
|H_\alpha(s,y)|
&\leq
a^{-1}\int_{\min(0,s)}^{\max(0,s)}
|\partial_{y_\alpha}\pi(t,y)|\,dt
+
|\partial_{y_\alpha}\psi(y)|
\\
&\leq
C{(R-r)}^{-M}\|h\|_{H^1(Q_R)}.
\end{aligned}
\]
Consequently,
\[
\|H\|_{L^\infty(Q_{\rho_6})}
\leq
C{(R-r)}^{-M}\|h\|_{H^1(Q_R)}.
\]

We next compute the first derivatives of $H$. For the $s$-component,
\[
\partial_sH_0
=
-a^{-1}\partial_s\pi(s,y),
\]
and
\[
\partial_{y_\beta}H_0
=
-a^{-1}
\partial_{y_\beta}\pi(s,y).
\]
For the transverse components, we have
\[
\partial_sH_\alpha
=
-a^{-1}\partial_{y_\alpha}\pi(s,y),
\]
and
\[
\partial_{y_\beta}H_\alpha
=
-a^{-1}\int_0^s
\partial_{y_\alpha y_\beta}\pi(t,y)\,dt
+
\partial_{y_\alpha y_\beta}\psi(y).
\]
It follows from the pressure and Poisson estimates that
\[
\|\nabla H\|_{L^\infty(Q_{\rho_6})}
\leq
C{(R-r)}^{-M}\|h\|_{H^1(Q_R)}.
\]

It remains to estimate the second derivatives of $H$. Differentiating
the formulas for the first derivatives, we obtain
\[
\partial_{ss}H_0
=
-a^{-1}\partial_{ss}\pi(s,y),
\]
\[
\partial_{s y_\beta}H_0
=
-a^{-1}\partial_{s y_\beta}\pi(s,y),
\]
and
\[
\partial_{y_\beta y_\gamma}H_0
=
-a^{-1}
\partial_{y_\beta y_\gamma}\pi(s,y).
\]
For the transverse components $H_\alpha$, we have
\[
\partial_{ss}H_\alpha
=
-a^{-1}\partial_{s y_\alpha}\pi(s,y),
\]
\[
\partial_{s y_\beta}H_\alpha
=
-a^{-1}\partial_{y_\alpha y_\beta}\pi(s,y),
\]
and
\[
\partial_{y_\beta y_\gamma}H_\alpha
=
-a^{-1}\int_0^s
\partial_{y_\alpha y_\beta y_\gamma}\pi(t,y)\,dt
+
\partial_{y_\alpha y_\beta y_\gamma}\psi(y).
\]

Every term in these formulas is controlled by the pressure and Poisson
estimates. In particular,
\[
\|D^2H\|_{L^\infty(Q_{\rho_6})}
\leq
C{(R-r)}^{-M}\|h\|_{H^1(Q_R)}.
\]
The factor $a^{-1}$ in front of the pressure terms cancels the factor
$a$ in the pressure bounds. Thus the estimate is uniform for
arbitrarily large values of the drift parameter $a$.

Since $Q_{\rho_6}$ has uniformly bounded measure, the estimates for
$H$, $\nabla H$, and $D^2H$ imply
\[
\begin{aligned}
\|H\|_{H^2(Q_{\rho_6})}
&=
\left(
\|H\|_{L^2(Q_{\rho_6})}^2
+
\|\nabla H\|_{L^2(Q_{\rho_6})}^2
+
\|D^2H\|_{L^2(Q_{\rho_6})}^2
\right)^{\frac12}
\\
&\leq
C\left(
\|H\|_{L^\infty(Q_{\rho_6})}
+
\|\nabla H\|_{L^\infty(Q_{\rho_6})}
+
\|D^2H\|_{L^\infty(Q_{\rho_6})}
\right)
\\
&\leq
C{(R-r)}^{-M}\|h\|_{H^1(Q_R)}.
\end{aligned}
\]
Therefore,
\[
\|H\|_{H^2(Q_{\rho_6})}
\leq
C{(R-r)}^{-M}\|h\|_{H^1(Q_R)},
\]
with a constant $C$ independent of $a\geq1$.
\end{proof}

\begin{lemma}\label{extension}
Let \(0<r<R\) be fixed.
There exists a linear extension operator \(\mathcal E\) such that,
for every smooth divergence-free vector field \(f\) on \(Q_R\),
\begin{align*}
\mathcal Ef=f\quad\text{in }Q_r,
\qquad
\mathcal Ef\in C_c^\infty(Q_R;\mathbb R^5),
\qquad
\nabla \cdot \mathcal Ef=0,
\end{align*}
and
\begin{align*}
\Y_a(\mathcal Ef)
\le C_{r,R}\Y_a(f;Q_R),
\qquad a\ge1.
\end{align*}
The constant \(C_{r,R}\) is independent of \(a\).
\end{lemma}

\begin{proof}
Choose \(\chi\in C_c^\infty(Q_R)\) such that
\begin{align*}
0\le\chi\le1,
\qquad
\chi=1\quad\text{in }Q_r.
\end{align*}
All derivatives of \(\chi\) used below have bounds depending
only on \(r,R\).

\textit{Step 1: construction of the skew-symmetric potential.}
Since \(Q_R\) is star-shaped with respect to the origin, define
\begin{align*}
A_{ij}(x)
=\int_0^1t^3
\bigl(x_jf_i(tx)-x_if_j(tx)\bigr)\,\mathrm dt,
\qquad i,j=1,\ldots,5.
\end{align*}
Then \(A_{ij}=-A_{ji}\). Moreover, using \(\nabla \cdot  f=0\),
\begin{align*}
\sum_{j=1}^5\partial_jA_{ij}(x)
&=\int_0^1t^3
\left(
5f_i(tx)+t\,x\cdot\nabla f_i(tx)
-f_i(tx)-t x_i\nabla \cdot  f(tx)
\right)\,\mathrm dt\\
&=\int_0^1
\left(
4t^3f_i(tx)+t^4x\cdot\nabla f_i(tx)
\right)\,\mathrm dt=\int_0^1
\frac{\mathrm d}{\mathrm dt}\bigl(t^4f_i(tx)\bigr)
\,\mathrm dt
=f_i(x).
\end{align*}
The term at \(t=0\) vanishes because \(f\) is smooth near
the origin.

Define
\begin{align*}
(\mathcal Ef)_i
=\sum_{j=1}^5\partial_j(\chi A_{ij})
=\chi f_i+\sum_{j=1}^5(\partial_j\chi)A_{ij},
\end{align*}
and extend this field by zero outside \(Q_R\).
It is smooth and compactly supported in \(Q_R\), and equals
\(f\) on \(Q_r\).
Skew symmetry gives
\begin{align*}
\begin{aligned}
\nabla \cdot \mathcal Ef
=\sum_{i,j=1}^5
\partial_i\partial_j(\chi A_{ij})=\frac12\sum_{i,j=1}^5
\partial_i\partial_j\bigl(\chi(A_{ij}+A_{ji})\bigr)
=0.
\end{aligned}
\end{align*}
The construction is linear in \(f\).

\textit{Step 2: slice estimates for \(A\) and \(\nabla A\).}
For \(0<t\le1\), \(|s|<R\), and \(p\in\{2,4\}\), the
transverse change of variables \(z=ty\) gives
\begin{align*}
\begin{aligned}
\|g(ts,t\,\cdot)\|_{L^p(B_R^4)}
&=\left(
\int_{B_R^4}|g(ts,ty)|^p\,\mathrm dy
\right)^{\frac1p}\\
&=t^{-\frac4p}
\left(
\int_{B_{tR}^4}|g(ts,z)|^p\,\mathrm dz
\right)^{\frac1p}\le t^{-\frac4p}
\|g\|_{L^\infty((-R,R);L^p(B_R^4))}.
\end{aligned}
\end{align*}
Since \(|x|\le\sqrt2R\) on \(Q_R\), Minkowski's integral
inequality yields
\begin{align*}
&\|A\|_{L^\infty((-R,R);L^4(B_R^4))}
\le C_R\int_0^1t^{3-1}\,\mathrm dt\,
\|f\|_{L^\infty((-R,R);L^4(B_R^4))}\le C_R
\|f\|_{L^\infty((-R,R);L^4(B_R^4))},\\
&\|A\|_{L^\infty((-R,R);L^2(B_R^4))}
\le C_R\int_0^1t^{3-2}\,\mathrm dt\,
\|f\|_{L^\infty((-R,R);L^2(B_R^4))}\le C_R
\|f\|_{L^\infty((-R,R);L^2(B_R^4))}.
\end{align*}

For each coordinate index \(\ell\),
\begin{align*}
\partial_\ell A_{ij}(x)
={}&\int_0^1t^3
\bigl(
\delta_{j\ell}f_i(tx)-\delta_{i\ell}f_j(tx)
\bigr)\,\mathrm dt+\int_0^1t^4
\bigl(
x_j\partial_\ell f_i(tx)
-x_i\partial_\ell f_j(tx)
\bigr)\,\mathrm dt.
\end{align*}
Applying the transverse \(L^2\) scaling to both integrals
gives
\begin{align*}
\|\nabla A\|_{L^\infty((-R,R);L^2(B_R^4))}
\le{}&
C\int_0^1t^{3-2}\,\mathrm dt\,
\|f\|_{L^\infty((-R,R);L^2(B_R^4))}\\
&+C_R\int_0^1t^{4-2}\,\mathrm dt\,
\|\nabla f\|_{L^\infty((-R,R);L^2(B_R^4))}\\
\le{}&
C_R\left(
\|f\|_{L^\infty((-R,R);L^2(B_R^4))}
+\|\nabla f\|_{L^\infty((-R,R);L^2(B_R^4))}
\right).
\end{align*}
The lower-order slice norm is controlled by
\begin{align*}
\|f\|_{L^\infty((-R,R);L^2(B_R^4))}
\le |B_R^4|^{\frac14}
\|f\|_{L^\infty((-R,R);L^4(B_R^4))}.
\end{align*}

Differentiating the expression for \(\mathcal Ef\), we have
\begin{align*}
\begin{aligned}
\partial_\ell(\mathcal Ef)_i
=\chi\,\partial_\ell f_i
+(\partial_\ell\chi)f_i+\sum_{j=1}^5(\partial_\ell\partial_j\chi)A_{ij}
+\sum_{j=1}^5(\partial_j\chi)\partial_\ell A_{ij}.
\end{aligned}
\end{align*}
Consequently,
\begin{align*}
\|\mathcal Ef\|_{L^\infty(\mathbb R;L^4(\mathbb R^4))}
&\le C_{r,R}\left(
\|f\|_{L^\infty((-R,R);L^4(B_R^4))}
+\|A\|_{L^\infty((-R,R);L^4(B_R^4))}
\right)\\
&\le C_{r,R}
\|f\|_{L^\infty((-R,R);L^4(B_R^4))},
\end{align*}
and
\begin{align*}
\|\nabla\mathcal Ef\|_{L^\infty(\mathbb R;L^2(\mathbb R^4))}
\le{}&
C_{r,R}\left(
\|\nabla f\|_{L^\infty((-R,R);L^2(B_R^4))}
+\|f\|_{L^\infty((-R,R);L^2(B_R^4))}
\right.\\
&\left.\hspace{14mm}
+\|\nabla A\|_{L^\infty((-R,R);L^2(B_R^4))}
+\|A\|_{L^\infty((-R,R);L^2(B_R^4))}
\right)\\
\le{}&
C_{r,R}\left(
\|f\|_{L^\infty((-R,R);L^4(B_R^4))}
+\|\nabla f\|_{L^\infty((-R,R);L^2(B_R^4))}
\right).
\end{align*}

\textit{Step 3: second derivatives of \(A\).}
Differentiating the formula for \(\partial_\ell A_{ij}\)
once more gives
\begin{align*}
\partial_m\partial_\ell A_{ij}(x)
={}&\int_0^1t^4
\bigl(
\delta_{j\ell}\partial_m f_i(tx)
-\delta_{i\ell}\partial_m f_j(tx)
\bigr)\,\mathrm dt\\ &+\int_0^1t^4
\bigl(
\delta_{jm}\partial_\ell f_i(tx)
-\delta_{im}\partial_\ell f_j(tx)
\bigr)\,\mathrm dt\\
&+\int_0^1t^5
\bigl(
x_j\partial_m\partial_\ell f_i(tx)
-x_i\partial_m\partial_\ell f_j(tx)
\bigr)\,\mathrm dt.
\end{align*}
Thus only first and second derivatives of \(f\) occur.

The full five-dimensional change of variables gives
\begin{align*}
\begin{aligned}
\|g(t\,\cdot)\|_{L^2(Q_R)}
&=\left(
\int_{Q_R}|g(tx)|^2\,\mathrm dx
\right)^{\frac12} =t^{-\frac52}
\left(
\int_{Q_{tR}}|g(z)|^2\,\mathrm dz
\right)^{\frac12}
\le t^{-\frac52}\|g\|_{L^2(Q_R)}.
\end{aligned}
\end{align*}
Hence
\begin{align*}
\|\nabla^2A\|_{L^2(Q_R)}
&\le
C\int_0^1t^{4-\frac52}\,\mathrm dt\,
\|\nabla f\|_{L^2(Q_R)}
+C_R\int_0^1t^{5-\frac52}\,\mathrm dt\,
\|\nabla^2f\|_{L^2(Q_R)}\\
&\le C_R\left(
\|\nabla f\|_{L^2(Q_R)}
+\|\nabla^2f\|_{L^2(Q_R)}
\right).
\end{align*}
Both integrals are finite, since
\begin{align*}
\int_0^1t^{\frac32}\,\mathrm dt=\frac25,
\qquad
\int_0^1t^{\frac52}\,\mathrm dt=\frac27.
\end{align*}
The lower-order full-space norms satisfy
\begin{align*}
\begin{aligned}
\|A\|_{L^2(Q_R)}+\|\nabla A\|_{L^2(Q_R)}
\le{}&(2R)^{\frac12}
\left(
\|A\|_{L^\infty((-R,R);L^2(B_R^4))}
+\|\nabla A\|_{L^\infty((-R,R);L^2(B_R^4))}
\right),\\
\|f\|_{L^2(Q_R)}
\le{}&(2R)^{\frac12}|B_R^4|^{\frac14}
\|f\|_{L^\infty((-R,R);L^4(B_R^4))},\\
\|\nabla f\|_{L^2(Q_R)}
\le{}&(2R)^{\frac12}
\|\nabla f\|_{L^\infty((-R,R);L^2(B_R^4))}.
\end{aligned}
\end{align*}

\textit{Step 4: the weighted second-derivative estimate.}
Using
\begin{align*}
(\mathcal Ef)_i
=\chi f_i+\sum_{j=1}^5(\partial_j\chi)A_{ij},
\end{align*}
the product rule and Steps 2--3 give
\begin{align*}
\|\nabla^2\mathcal Ef\|_{L^2(\mathbb R^5)}
\le C_{r,R}
\bigl(\|f\|_{H^2(Q_R)}+\|A\|_{H^2(Q_R)}\bigr).
\end{align*}
Hence, since \(a\ge1\),
\begin{align*}
a^{-\frac12}\|\nabla^2\mathcal Ef\|_{L^2(\mathbb R^5)}
\le C_{r,R}\bigl(
&a^{-\frac12}\|\nabla^2f\|_{L^2(Q_R)}
+\|f\|_{L^\infty((-R,R);L^4(B_R^4))}\\
&+\|\nabla f\|_{L^\infty((-R,R);L^2(B_R^4))}
\bigr).
\end{align*}
Combining this with the slice estimates proves
\(\Y_a(\mathcal Ef)\le C_{r,R}\Y_a(f;Q_R)\).
\end{proof}

\end{document}